\documentclass[11 pt,a4paper,leqno]{article}
\usepackage{latexsym,amsmath,amscd,amsfonts}
\usepackage[francais,english]{babel}

\RequirePackage{amsmath}
\RequirePackage{amssymb}
\RequirePackage{ifthen}
\RequirePackage{theorem}
\RequirePackage{pslatex}

\title{\sc \'Equations aux $q$-diff\'erences et fibr\'es vectoriels 
holomorphes sur la courbe elliptique $\C^{*}/q^{\Z}$}

\author{Jacques Sauloy
\footnote{Laboratoire Emile Picard, CNRS UMR 5580, U.F.R. M.I.G.,
118, route de Narbonne, 31062 Toulouse CEDEX 4}}

\date{}

\def\sq{\sigma_q}

\def\C{{\mathbf C}}
\def\Q{{\mathbf Q}}
\def\Z{{\mathbf Z}}
\def\R{{\mathbf R}}
\def\N{{\mathbf N}}

\def\Sr{{\mathbf{S}}}

\def\L{\mathcal{L}}
\def\ii{{\text{i}}}

\def\F{{\mathcal{F}}}
\def\P{{\mathcal{P}}}
\def\B{{\mathcal{B}}}
\def\Bd{{\mathcal{B}_{q}^{(\delta)}}}

\def\O{{\mathcal{O}}}
\def\M{{\mathcal{M}}}
\def\Ma{{\text{Mat}}}
\def\Co{{\mathcal{C}}}
\def\D{{\mathcal{D}_{q,K}}}
\def\gr{{\text{gr}}}

\def\Eq{{\mathbf{E}_{q}}}
\def\G{{\mathfrak{G}}}
\def\g{{\mathfrak{g}}}
\def\V{{\mathfrak{V}}}
\def\X{{\mathcal{X}}}

\def\Kr{{\C(z)}}

\def\Ka{{\C(\{z\})}}

\def\Raq(d){{\C\{\xi\}_{q,(\delta)}}}
\def\Kaq(d){{\C(\{\xi\})_{q,(\delta)}}}

\def\Kf{{\C((z))}}
\def\Rw{{\mathcal{O}(\C^{*})}}
\def\Kw{{\mathcal{M}(\C^{*})}}

\def\DM{{DiffMod\left(K,\sq\right)}}
\def\DMr{{DiffMod\left(\Kr,\sq\right)}}
\def\DMa{{DiffMod\left(\Ka,\sq\right)}}
\def\DMf{{DiffMod\left(\Kf,\sq\right)}}

\def\Er{{{\mathcal{E}}}}

\def\E{{{\mathcal{E}}^{(0)}}}
\def\Ef{{{\mathcal{E}}^{(0)}_{f}}}
\def\Ee{{{\mathcal{E}}^{(0)}_{1}}}
\def\Epe{{{\mathcal{E}}^{(0)}_{p,1}}}

\def\Der{{\dot{\Delta}}}
\def\Derc{{\Der_{\overline{c}}}}

\def\cad{{c'est-\`a-dire}}
\def\ie{{\emph{i.e.}}}
\def\cf{{\emph{cf.}}}

\def\div{\text{div}}

\def\Im{\text{Im}}

\def\Aut{\text{Aut}}
\def\GL{\mathcal{GL}}
\def\irr{\text{irr}}
\def\Sp{\text{Sp}}
\def\cl{\text{cl}}
\def\Res{\text{Res}}

\def\St{\mathfrak{St}}
\def\st{\mathfrak{st}}

\newtheorem{thm}{Th\'eor\`eme}[section]
\newtheorem{lemma}[thm]{Lemme}
\newtheorem{prop}[thm]{Proposition}
\newtheorem{cor}[thm]{Corollaire}

\def\Pr{\textsl{D\'emonstration. - }}
\def\Ex{\noindent\textbf{Exemple.~}}
\def\Exs{\noindent\textbf{Exemples. \\}}
\def\Rem{\noindent\textbf{Remarque.~}}
\def\Rems{\noindent\textbf{Remarques. \\}}

\begin{document}

\maketitle

\bigskip \hrule \bigskip

\centerline{\textbf{\emph{R\'esum\'e}}}

\emph{Nous pr\'esentons diverses applications des fibr\'es vectoriels
aux \'equations aux $q$-diff\'erences, dans la lign\'ee de la
correspondance de Weil.}

\bigskip \hrule \bigskip

\centerline{\textbf{\emph{Abstract}}}

\emph{We present some applications of vector bundles to $q$-difference
equations, in continuation of Weil's correspondance.}

\bigskip \hrule \bigskip

\tableofcontents

\bigskip \hrule 



\section{Introduction}
\label{section:introduction}

Divers fils math\'ematiques et historiques relient les \'equations
aux $q$-diff\'erences aux \emph{fibr\'es vectoriels holomorphes sur
une courbe elliptique}
\footnote{Dans tout le texte, nous dirons ``fibr\'e'' pour ``fibr\'e 
vectoriel holomorphe'' (sur une surface de Riemann).}. 
Ces derni\`eres ann\'ees, ces derniers sont apparus \`a plusieurs 
reprises comme un cadre naturel pour des probl\`emes de classification 
et de th\'eorie de Galois (probl\`eme de Riemann-Hilbert). Il est peut-\^etre
temps de survoler et de mettre en ordre des r\'esultats \'epars,
dont certains ont \'et\'e \'enonc\'es dans diverses conf\'erences
(Groningen, Conf\'erence Ramis, Lisbonne, Luminy, Kyoto, Tordesillas)
mais n'ont jamais \'et\'e publi\'es. Ces r\'esultats ont \'et\'e
tr\`es largement motiv\'es par les travaux de Ramis, Zhang et l'auteur
et l'une des raisons de non publication est le blocage sur une question
difficile, celle du ``probl\`eme global'' 
(section \ref{section:constructionsglobales}). Cependant les perc\'ees
des derni\`eres ann\'ees sur le probl\`eme local (\cite{RSZ},
\cite{RS1} et \cite{RS2}) nous encouragent. \\

L'article comprend peu de r\'esultats extraordinaires mais permet 
un \'eclairage nouveau de la th\'eorie. Il permet en particulier
de proposer une \'enigme (apparition de la dualit\'e de Serre)
et un probl\`eme ouvert (le probl\`eme global mentionn\'e ci-dessus). 
Nous n'\'evoquons pas deux autres pistes, celle de la \emph{confluence} 
(\cite{JSAIF}, \cite{JSGAL}) et celle des \emph{d\'eformations 
isomonodromiques} (\cite{JSISO}). \\

Nous nous occupons principalement d'\'equations aux $q$-diff\'erences
et ne sommes venus aux fibr\'es vectoriels que par n\'ecessit\'e:
nous ne pr\'etendons \`a aucune expertise dans ce domaine, et esp\'erons
au contraire que les sp\'ecialistes nous apporteront leurs lumi\`eres. \\

{\small
Ce fut un plaisir tout particulier de parler de tout cela \`a la
conf\'erence en l'honneur de Jose-Manuel Aroca, Gran Jefe Capitan
Pirata, en pr\'esence de tant d'amis de Valladolid et d'ailleurs.
\`A Valladolid et \`a Tordesillas, on rit beaucoup avant, pendant
et apr\`es les expos\'es (parfois, \`a la place) parce que le plaisir 
de faire des math\'ematiques s'y exprime plus librement qu'ailleurs. 
Merci pour tout cela \`a Jose-Manuel, l'\^ame du groupe. \\
J'avais pr\'efac\'e mon expos\'e (en anglais) \`a Tordesillas de 
la d\'edicace suivante: \\
\emph{With a special thought for Jean Giraud,} \\
\emph{who, a long time ago, guided my first steps} \\
\emph{into the wild world of singularities ...} \\
Jean Giraud, qui n'avait pu assister \`a la conf\'erence, nous a quitt\'es 
le 27 mars 2007. Je partage ici ma tristesse avec nos amis espagnols.
}


\subsection{Apparition des fibr\'es dans 
la th\'eorie des \'equations fonctionnelles}
\label{subsection:apparition}

Le th\'eor\`eme cl\'e dans la r\'esolution par Birkhoff du probl\`eme 
de Riemann-Hilbert (\cite{Birkhoff1}) est un th\'eor\`eme de factorisation 
de matrice holomorphe. Dans \cite{Rohrl57}, \cite{Rohrl62}, R\"ohrl 
a interpr\'et\'e 
ce th\'eor\`eme en termes de \emph{trivialit\'e de fibr\'e vectoriel}
(voir aussi \cite{Forster}). Dans \cite{vdPS}, van der Put et Singer
donnent de cette factorisation une preuve moderne, qui s'appuie
directement sur la cohomologie des fibr\'es vectoriels sur une surface
de Riemann, et l'appliquent (dans la droite ligne de \cite{Birkhoff1})
aux \'equations aux diff\'erences et aux $q$-diff\'erences. Auparavant,
Praagman, un \'el\`eve de van der Put, avait invoqu\'e la trivialit\'e
m\'eromorphe des fibr\'es pour d\'emontrer l'existence d'un syst\`eme
fondamental de solutions m\'eromorphes sur $\C^{*}$ pour les \'equations
aux diff\'erences et aux $q$-diff\'erences (\cite{Praagman}). Cependant,
dans tous ces cas, les fibr\'es n'interviennent qu'\`a travers leurs
propri\'et\'es cohomologiques, et non en tant qu'objets g\'eom\'etriques. \\

Dans \cite{BG}, Baranovsky et Ginzburg \'etudient la classification
formelle des \'equations aux $q$-diff\'erences fuchsiennes (dans une 
autre terminologie, li\'ee aux groupes de lacets). Ils caract\'erisent
chaque classe \emph{formelle} \`a l'aide d'un objet \emph{analytique},
un fibr\'e vectoriel sur la courbe elliptique $\Eq = \C^{*}/q^{\Z}$.
Sur une suggestion de Kontsevitch, ils en d\'eduisent le groupe de
Galois local. Ind\'ependamment, l'auteur a obtenu dans \cite{JSGAL}
la classification (formelle ou analytique, ce qui revient au m\^eme
dans ce cas) des \'equations aux $q$-diff\'erences fuchsiennes par
des fibr\'es plats, d'o\`u se d\'eduit la description compl\`ete
du groupe de Galois local et celle moins d\'etaill\'ee du groupe
de Galois global (cas ab\'elien r\'egulier). \\

Nous allons, dans cette introduction, suivre le chemin inverse 
et montrer comment la description des fibr\'es sur une courbe 
(resp. une courbe elliptique) se traduit naturellement en termes 
d'\'equations fonctionnelles (resp. d'\'equations aux $q$-diff\'erences).


\subsubsection{La correspondance de Weil}
\label{subsubsection:correspondancedeWeil}

Dans \cite{Weil}, Weil propose, sous le nom de $G$-diviseurs, 
une g\'en\'eralisation non-ab\'elienne de la notion de diviseur
sur une surface de Riemann. Ces $G$-diviseurs ne sont autres 
que des fibr\'es vectoriels avant la lettre. Selon la pr\'esentation
``moderne'' de \cite{Grothendieck} (et sous une forme simplifi\'ee),
cela donne ce qui suit.

\paragraph{Fibr\'es \'equivariants.}

Soit $E$ une surface de Riemann, et soit $\tilde{E}$ son  
rev\^etement universel, qui est donc \'egalement une surface
de Riemann. Nous noterons $\pi: \tilde{E} \rightarrow E$ la
projection canonique. \\

Soit $\F$ un fibr\'e (vectoriel holomorphe) sur $E$. En relevant $\F$ 
\`a $\tilde{E}$ par $\pi$, on obtient un fibr\'e $\tilde{\F} = \pi^{*} \F$, 
qui est trivial puisque $\tilde{E}$ est simplement connexe. 
On \'ecrit donc $\tilde{\F} = \tilde{E} \times V$, o\`u $V$ est
un $\C$-espace vectoriel de dimension finie. Provenant de $E$, 
ce fibr\'e trivial est muni d'une \emph{action \'equivariante} du groupe
$G = \Aut(\tilde{E} / E) = \pi_{1}(E)$ (nous ne pr\'ecisons
pas le point-base pour le groupe fondamental $\pi_{1}$,
qui n'apparaitra qu'en tant que groupe des automorphismes
du rev\^etement). Le mot ``action \'equivariante'' signifie ici
``action sur $\tilde{E} \times V$ qui commute avec l'action
sur $\tilde{E}$'' (on dit aussi que $\tilde{\F}$ est un fibr\'e
\'equivariant). Une telle action est compl\`etement d\'ecrite
par l'action naturelle $(\gamma,x) \mapsto \gamma.x$ de $G$ sur $\tilde{E}$ 
et par la donn\'ee d'une application holomorphe (en la seconde variable):
$$
A: G \times \tilde{E} \longrightarrow \GL(V).
$$
Tout $\gamma \in G$ op\`ere alors sur $\tilde{\F} = \tilde{E} \times V$
par l'application:
$$
(x,X) \mapsto \left(\gamma.x,A(\gamma,x) X\right).
$$
Pour que ce soit bien une op\'eration de groupe, il faut,
et il suffit, que soit r\'ealis\'ee une condition de cocycle:
$$
\forall \gamma,\gamma' \in G \;,\; \forall x \in \tilde{E} \;,\;
A(\gamma' \gamma,x) = A(\gamma',\gamma.x) A(\gamma,x).
$$
On peut \'egalement exprimer, par une condition de cobord, la
trivialit\'e du fibr\'e $\F$ de d\'epart ou, plus g\'en\'eralement,
\`a quelle condition deux cocycles repr\'esentent des fibr\'es isomorphes. \\

Un morphisme $\F \rightarrow \F'$ de fibr\'es sur $E$ se rel\`eve
en un morphisme 
$\tilde{\F} = \tilde{E} \times V \rightarrow 
\tilde{\F'} = \tilde{E} \times V'$ 
de fibr\'es sur $\tilde{E}$ compatible avec la structure ci-dessus:
si $\tilde{\F}$ et $\tilde{\F'}$ sont respectivement d\'ecrits par
les cocycles $A$ et $A'$, le morphisme $\tilde{\F} \rightarrow \tilde{\F'}$ 
est de la forme $(x,X) \mapsto \left(x,F(x) X\right)$, o\`u $F$ est
une application holomorphe de $\tilde{E}$ dans $\L(V,V')$, qui satisfait
\`a la condition suivante:
$$
\forall \gamma \in G \;,\; \forall x \in \tilde{E} \;,\;
F(\gamma.x) A(\gamma,x) = A'(\gamma,x) F(x).
$$

\paragraph{Description g\'eom\'etrique.}

Supposons r\'eciproquement donn\'e le cocycle holomorphe  
(en la seconde variable)
$A: \pi_{1}(E) \times \tilde{E} \rightarrow \GL(V)$. On lui associe 
la relation d'\'equivalence $\sim_{A}$ sur le fibr\'e trivial 
$\tilde{F} = \tilde{E} \times V$ engendr\'ee par les relations:
$(x,X) \sim_{A} \left(\gamma.x,A(\gamma,x) X\right)$: la relation
$\sim_{A}$ provient donc d'une action \'equivariante de $\pi_{1}(E)$ 
sur $\tilde{E}$. En un sens \'evident, cette relation est compatible
avec la relation $\sim$ sur $\tilde{E}$ induite par l'action 
de $\pi_{1}(E)$. Le fibr\'e sur $E$ associ\'e, que nous noterons 
$\tilde{F}_{A}$, s'obtient par passage au quotient de la projection
$\tilde{F} = \tilde{E} \times V \rightarrow \tilde{E}$ par ces relations
d'\'equivalence:
$$
\F_{A} = \dfrac{\tilde{E} \times V}{\sim_{A}} \longrightarrow
E = \dfrac{\tilde{E}}{\sim} \cdot
$$
On peut alors d\'ecrire le faisceau des sections de $\F_{A}$.
Soit $V$ un ouvert de $E$. Alors l'espace des sections de $\F_{A}$
sur $V$ est:
$$
\Gamma(V,\F_{A}) = 
\{X: \pi^{-1}(V) \rightarrow V \text{~holomorphes~} \mid
\forall x \in \pi^{-1}(V) \;,\; \forall \gamma \in \pi_{1}(E) \;,\;
X(\gamma.x) = A(\gamma,x) X(x) \}.
$$

\Ex
Prenons $E = \C^{*}$. Alors $\tilde{E} = \C$ sur lequel 
$\pi_{1}(E) = \Z$ agit par translations, et la projection canonique
$\pi: \tilde{E} \rightarrow E$ est ici $x \mapsto e^{2 \ii \pi x}$.
La condition de cocycle entraine que $A$ est enti\`erement
d\'etermin\'ee par la matrice $A(1,x)$. Notons (abusivement)
$A(x) = A(1,x)$. De m\^eme, la condition qui d\'efinit les
sections peut se tester simplement en prenant $\gamma = 1$:
$$
\Gamma(V,\F_{A}) = 
\{X: \pi^{-1}(V) \rightarrow V \text{~holomorphes~} \mid
\forall x \in \pi^{-1}(V) \;,\; X(x+1) = A(x) X(x)\}.
$$
On voit bien la parent\'e avec les \'equations fonctionnelles. \\

Si l'on note $\underline{1} = \F_{1}$ (``objet unit\'e'') le fibr\'e
en droites trivial sur $E$, associ\'e au cocycle trivial
$(\gamma,x) \mapsto 1 \in \GL(\C)$, le lecteur pourra v\'erifier
que les morphismes de $\underline{1}$ dans un fibr\'e $\F_{A}$
quelconque s'identifient aux sections globales de $\F_{A}$.

\paragraph{Fibr\'es plats et repr\'esentations de $\pi_{1}(E)$.}

Un cas important est celui o\`u, \`a isomorphisme pr\`es, on peut
supposer $A(\gamma,x)$ ind\'ependant de $x \in \tilde{E}$: on
l'\'ecrit donc $A(\gamma)$, et la condition de cocycle dit alors
que $\gamma \mapsto A(\gamma)$ est une repr\'esentation de
$\pi_{1}(E)$ dans $\GL(V)$. Un tel fibr\'e est appel\'e \emph{plat}
(\cite{Gunning}). Les fibr\'es plats admettent une caract\'erisation
topologique: les classes de Chern sur leurs facteurs ind\'ecomposables
sont nulles; et une caract\'erisation diff\'erentielle: on peut les
munir d'une connexion holomorphe. \emph{Nous n'aurons pas l'usage
de ces caract\'erisations}
\footnote{Van der Put et Reversat utilisent la seconde dans \cite{vdPR},
voir l\`a dessus la section \ref{subsection:equationsirregulieres}.}.
On obtient ainsi la c\'el\`ebre \emph{correspondance de Weil} entre
fibr\'es plats et repr\'esentations du groupe fondamental. \\

Il faut cependant prendre garde que cette correspondance n'est pas
une \'equivalence entre la cat\'egorie des fibr\'es plats sur $E$
et celle des repr\'esentations de $\pi_{1}(E)$. Soient en effet
$A: \pi_{1}(E) \rightarrow \GL(V)$ et $A': \pi_{1}(E) \rightarrow \GL(V')$
deux telles repr\'esentations, et soient $\F_{A}$ et $\F_{A'}$ 
les fibr\'es plats qui leur correspondent respectivement. Un
morphisme de $\F_{A}$ dans $\F_{A'}$ est d\'ecrit comme une
application holomorphe $F: \tilde{E} \rightarrow \L(V,V')$, 
telle que:
$$
\forall \gamma \in G \;,\; \forall x \in \tilde{E} \;,\;
F(\gamma.x) A(\gamma) = A'(\gamma) F(x).
$$
Si $F$ est constant sur $\tilde{E}$, c'est bien un morphisme 
de repr\'esentations, mais pas autrement. Nous en verrons un exemple 
\`a la section suivante, et des cons\'equences pour le groupe de Galois 
\`a la section \ref{subsubsection:equationsfuchsiennes}.


\subsubsection{Le cas des fibr\'es sur une courbe elliptique}
\label{subsubsection:fibressurEq}

\paragraph{Fibr\'es sur $\C/(\Z + \Z \tau)$.}

Prenons pour $E$ la courbe elliptique 
\footnote{\emph{A priori}, la surface de riemann $E$ devrait \^etre
appel\'ee ``tore complexe'', mais l'on sait que c'est essentiellement
la m\^eme chose qu'une courbe elliptique.}
$\C/\Lambda_{\tau}$, o\`u
$\Im \, \tau < 0$ et $\Lambda_{\tau} = \Z + \Z \tau$. (Nous poserons
plus loin $q = e^{2 \ii \pi \tau}$ et voudrons avoir $| q | > 1$.)
Ici, $\pi_{1}(E) = \Lambda_{\tau}$ agit sur $\tilde{E} = \C$ par
translations. Notons encore $\pi: \C \rightarrow E$ la projection
canonique. Pour tout cocycle $A$, notons $A_{1}(x) = A(1,x)$
et $A_{\tau}(x) = A(\tau,x)$. \`A cause de la relation de commutation
$\tau + 1 = 1 + \tau$, la condition de cocycle entraine:
$$
\forall x \in \C \;,\; A_{\tau}(x+1) A_{1}(x) = A_{1}(x+\tau) A_{\tau}(x).
$$
R\'eciproquement, deux applications holomorphes de $\C$ dans $\GL(V)$
qui v\'erifient cette relation s'\'etendent de mani\`ere unique en
un cocycle $A$ et d\'efinissent donc un fibr\'e $\F = \F_{A}$ sur $E$.
Les sections de ce fibr\'e sur l'ouvert $V \subset E$ s'identifient 
aux solutions holomorphes sur $\pi^{-1}(V) \subset \C$ de l'\'equation 
fonctionnelle:
$$
\forall x \in \pi^{-1}(V) \;,\;
X(x+1) = A_{1}(x) X(x) \quad \text{et} \quad X(x+\tau) = A_{\tau}(x) X(x).
$$
Si $\F' = \F_{A'}$ est le fibr\'e d\'efini par $A'_{1}$ et $A'_{\tau}$
(holomorphes de $\C$ dans $\GL(V')$), un morphisme de $\F$ dans $\F'$
est repr\'esent\'e par  une application holomorphe de $\C$ dans
$\L(V,V')$ telle que:
$$
\forall x \in \C \;,\;
F(x+1) A_{1}(x) = A'_{1}(x) F(x) \quad \text{et} \quad 
F(x+\tau) A'_{\tau}(x) = A_{\tau}(x) F(x).
$$
Le fibr\'e $\F_{A}$ est plat si, \`a isomorphisme pr\`es,
on peut supposer que $A_{1}$ et $A_{\tau}$ ne d\'ependent pas 
de $x$: $A_{1},A_{\tau} \in \GL(V)$. La condition de cocycle
dit alors que ces deux matrices commutent. La repr\'esentation
de $\pi_{1}(E) = \Lambda_{\tau}$ associ\'ee \`a $\F_{A}$ par
la correspondance de Weil est celle d\'efinie par $1 \mapsto A_{1}$
et $\tau \mapsto A_{\tau}$.

\paragraph{Fibr\'es sur $\C^{*}/q^{\Z}$.}
\label{paragraph:fibressurEq}

Pour trivialiser le fibr\'e $\F$ sur $E$, il n'est cependant pas
n\'ecessaire de le relever au rev\^etement universel $\C$. Ce
rev\^etement se factorise en 
$\C \rightarrow \C/\Z \rightarrow \C/\Lambda_{\tau}$. Or, l'application
$x \mapsto z = e^{2 \ii \pi x}$ permet d'identifier $\C/\Z$ \`a la
surface de Riemann ouverte $\C^{*}$. La m\^eme application permet
d'identifier $E = \C/\Lambda_{\tau}$ \`a $\Eq = \C^{*}/q^{\Z}$, o\`u
$q = e^{2 \ii \pi \tau}$ est un nombre complexe arbitraire de module
$|q| > 1$. On peut alors relever le fibr\'e $\F$ sur $\Eq$ en un fibr\'e
sur $\C^{*}$ par le rev\^etement $\C^{*} \rightarrow \Eq$. L'inter\^et 
de cette op\'eration est que tout fibr\'e vectoriel holomorphe sur une
surface de Riemann ouverte (\ie\ non compacte) est trivial
(\cite{Gue-Nar}, th\'eor\`eme 3 p. 184). \\

Le formalisme des fibr\'es \'equivariants d\'ecrit \`a la section
\ref{subsubsection:correspondancedeWeil} s'applique alors tout aussi 
bien ici. Nous partirons donc maintenant de la description ``de Jacobi'' 
(ou ``de Tate'') des courbes elliptiques pour fixer nos notations.
Soit $q$ un complexe de module $|q| > 1$. Soit $\Eq = \C^{*}/q^{\Z}$.
On note $\pi$ la projection canonique $\C^{*} \rightarrow \Eq$.
C'est un rev\^etement, dont le groupe $\Aut(\C^{*}/\Eq)$ est $q^{\Z}$
agissant sur $\C^{*}$ en tant que sous-groupe. \\

Tout fibr\'e $\F$ sur $\Eq$ se rel\`eve par $\pi$ en un fibr\'e trivial 
$\tilde{\F} = \C^{*} \times V$ muni d'une action \'equivariante de $q^{\Z}$,
autrement dit, d'une application holomorphe (en la seconde variable)
$A: q^{\Z} \times \C^{*} \rightarrow \GL(V)$. Celle-ci satisfait la
condition de cocycle suivante:
$$
\forall m,n \in \Z \;,\; \forall z \in \C^{*} \;,\;
A(q^{m+n},z) = A(q^{m},q^{n} z) A(q^{n},z).
$$
Il est ais\'e de voir que la donn\'ee de $A(q,z)$ d\'etermine $A$.
Notant abusivement $A(z) = A(q,z)$, on trouve que l'on a, pour 
$n \geq 1$: $A(q^{n},z) = A(q^{n-1} z) \cdots A(z)$; et, pour 
$n \leq -1$ ... une formule laiss\'ee en exercice au lecteur !
Ainsi, il revient au m\^eme de se donner un cocycle $A$ ou une 
application holomorphe $A: \C^{*} \rightarrow \GL(V)$; et une telle
fonction matricielle $A$ d\'efinit un fibr\'e $\F_{A}$ sur $\Eq$.
Ce dernier peut \^etre construit g\'eom\'etriquement ainsi:
$$
\F_{A} = \dfrac{\C^{*} \times V}{\sim_{A}} \longrightarrow
\Eq = \dfrac{\C^{*}}{\sim},
$$
o\`u les relations d'\'equivalences sont d\'efinies par
$(z,X) \sim_{A} \left(q z,A(z) X\right)$ et $z \sim q z$.
Une section de $\F_{A}$ sur l'ouvert $V \subset \Eq$
s'identifie \`a une solution holomorphe sur $\pi^{-1}(V)$
de l'\emph{\'equation aux $q$-diff\'erences}:
$$
X(q z) = A(z) X(z).
$$
Soient $A: \C^{*} \rightarrow \GL(V)$, $A': \C^{*} \rightarrow \GL(V')$
deux telles applications holomorphes et $\F = \F_{A}$, $\F' = \F_{A'}$
les fibr\'es sur $\Eq$ associ\'es. Un morphisme de $\F$ dans $\F'$ 
est repr\'esent\'e par une application holomorphe
$F : \C^{*} \rightarrow \L(V,V')$ telle que:
$$
\forall z \in \C^{*} \;,\; F(q z) A(z) = A'(z) F(z).
$$
Par exemple, si l'on note $\underline{1}$ le fibr\'e en droites
trivial
\footnote{La notation $\underline{1}$ d\'esigne l'objet unit\'e,
c'est \`a dire le neutre pour le produit tensoriel, dans une
``cat\'egorie tannakienne''.}, 
provenant de la fonction constante $1$ de $\C^{*}$ dans
$\C^{*} = \GL(\C)$, on voit que les morphismes de $\underline{1}$
dans $\F$ s'identifient aux sections de $\F$. \\

\Rem
Si l'on rel\`eve le fibr\'e $\F$ sur $\Eq$ d'abord \`a $\C^{*}$
puis \`a $\C$, on obtient successivement $\C^{*} \times V$,
(muni d'une fonction matricielle $A(z)$ sur $\C^{*}$), et $\C \times V$. 
Ainsi, le fibr\'e trivial \'equivariant sur $\C$ d\'ecrit plus haut \`a 
l'aide des fonctions matricielles $A_{1}$ et $A_{\tau}$ sur $\C$,
peut-il toujours \^etre r\'ealis\'e en prenant $A_{1}(x) = Id_{V}$
et $A_{\tau}(x) = A(e^{2 \ii \pi x})$. Pour \^etre pr\'ecis, parmi
toutes les trivialisations de l'image r\'eciproque de $\F$ sur $\C$,
l'une au moins est munie d'une action \'equivariante de cette nature.
Cette propri\'et\'e, qui traduit la trivialit\'e des fibr\'es holomorphes
sur $\C^{*}$, \'equivaut \`a la suivante: la fonction matricielle
$A_{1}$ \'etant donn\'ee, l'\'equation fonctionnelle 
$X(x+1) = A_{1}(x) X(x)$ admet une solution fondamentale (\cad\
une solution \`a valeurs dans $\GL(V)$) holomorphe. \\
Si l'on se restreint aux fibr\'es plats, on en d\'eduit (correspondance
de Weil) que toute repr\'esentation de $\Z^{2} \simeq \Z + \Z \tau$ est 
\'equivalente \`a une repr\'esentation triviale sur le premier facteur.
C'est \'evidemment faux pour l'\'equivalence habituelle des repr\'esentations,
mais c'est vrai au sens de l'\'equivalence ``\'equivariante'' d\'ecrite
\`a la fin de la section \ref{subsubsection:correspondancedeWeil}.

\paragraph{Relations avec la th\'eorie classique 
des \'equations fonctionnelles.}

Comme on l'a vu, les sections de $\F_{A}$ s'identifient aux solutions 
de l'\emph{\'equation aux $q$-diff\'erences}: $X(qz) = A(z) X(z)$.
Il y a cependant une diff\'erence notable avec la th\'eorie classique
des \'equations fonctionnelles (\cite{Birkhoff1}, \cite{DRSZ},
\cite{GR}, \cite{vdPS}, \cite{JSAIF}): ici, la matrice $A(z)$ est
holomorphe sur $\C^{*}$, et m\^eme r\'eguli\`ere, \ie\ 
son inverse $A^{-1}$ est aussi holomorphe; alors que dans la th\'eorie 
classique, la matrice $A(z)$ est rationnelle (et inversible). Ainsi:
\begin{itemize}
\item{Pour ramener la th\'eorie classique \`a celle des fibr\'es,
il faut se d\'ebarrasser des p\^oles de $A$ et de $A^{-1}$.}
\item{Pour ramener la th\'eorie des fibr\'es sur $\Eq$ \`a 
la th\'eorie classique des \'equations aux $q$-diff\'erences,
il faut dompter la sauvagerie des \'equations (et des solutions)
en $0$ et en $\infty$.}
\end{itemize}
Comme on le verra (section \ref{subsubsection:equationsfuchsiennes}),
la th\'eorie fuchsienne vient naturellement se placer \`a l'intersection
des deux points de vue. 

\paragraph{Le cas des \'equations aux diff\'erences.}

Dans le cas des \'equations aux $q$-diff\'erences, le corps 
des constantes de la th\'eorie (solutions m\'eromorphes sur $\C^{*}$ 
de l'\'equation triviale $f(qz) = f(z)$) s'identifie au corps 
des fonctions elliptiques $\M(\Eq)$, corps des fonctions m\'eromorphes 
sur la surface de Riemann compacte $\Eq$; celle-ci s'identifie \`a 
une courbe alg\'ebrique (courbe elliptique) et $\M(\Eq)$ \`a un corps 
de fonctions alg\'ebriques; plus g\'en\'eralement, les fibr\'es vectoriels
holomorphes sont alg\'ebriques (\cite{SerreGAGA}). \\
La th\'eorie des \emph{\'equations aux diff\'erences} $X(z+1) = A(z) X(z)$ 
se pr\^ete \'egalement au point de vue des fibr\'es, mais c'est plus
compliqu\'e. En effet, la surface de Riemann appropri\'ee est ici
$E = \C/\Z \simeq \C^{*}$, mais celle-ci n'est pas compacte. 
Le corps des constantes de la th\'eorie (solutions m\'eromorphes sur $\C$ 
de l'\'equation triviale $f(z+1) = f(z)$) est ``tr\`es gros''.
Il faut donc artificiellement imposer des conditions de croissance
aux solutions pour les maitriser. Au fond, le cas des \'equations
aux diff\'erences est une d\'eg\'enerescence du cas des \'equations
aux $q$-diff\'erences. C'est parce que l'op\'erateur de translation
$z \mapsto z+1$ n'a qu'un point fixe sur la sph\`ere de Riemann,
alors que l'op\'erateur de dilatation $z \mapsto q z$ en a deux.
Anne Duval (\cite{AD}, voir aussi \cite{ADJR}) a \'etudi\'e la
\emph{confluence} de ces deux points fixes en un seul et ses 
cons\'equences sur les liens entre les deux types d'\'equations.


\subsection{Conventions g\'en\'erales}
\label{subsection:conventions}

Dans tout l'article, nous fixerons un nombre complexe $q \in \C$
de module $| q | > 1$. Nous noterons $\Eq$ la courbe elliptique
$\C^{*}/q^{\Z}$ et $\pi: \C^{*} \rightarrow \Eq$ la projection
canonique. L'image dans $\Eq$ de $a \in \C^{*}$ sera not\'ee
$\overline{a}$. La \emph{spirale logarithmique discr\`ete}
$\pi^{-1}(\overline{a}) = a q^{\Z}$ sera not\'ee $[a;q]$. 
On \'ecrira alors $[a,b;q] = [a;q] \cup [b;q]$, etc. \\

L'op\'erateur de dilatation $z \mapsto q z$ de la sph\`ere 
de Riemann $\Sr$ induit un automorphisme $\sq$ sur de nombreux 
anneaux ou corps de fonctions, par la formule $(\sq f)(z) = f(qz)$ 
(cette notation s'\'etend naturellement \`a des vecteurs ou 
des matrices de fonctions). Les principaux corps d'inter\^et 
sont $\Kr$ (fonctions rationnelles), $\Ka$ (germes m\'eromorphes 
en $0$), $\Kf$ (s\'eries de Laurent formelles) et $\Kw$ (fonctions 
m\'eromorphes sur $\C^{*}$). Plus g\'en\'eralement, le corps des fonctions 
m\'eromorphes (resp. l'anneau des fonctions holomorphes) sur une surface 
de Riemann $E$ est not\'e $\M(E)$ (resp. $\O(E)$).


\subsubsection{Fonctions}
\label{subsubsection:fonctions}

Les fonctions m\'eromorphes sur $E = \C/(\Z + \Z \tau)$ s'identifient
aux fonctions m\'eromorphes sur $\C$ admettant le r\'eseau de
p\'eriodes $\Z + \Z \tau$: c'est la description classique du corps
$\M(E)$ des fonctions elliptiques. Les fonctions m\'eromorphes sur 
$\Eq = \C^{*}/q^{\Z}$ s'identifient de m\^eme aux fonctions m\'eromorphes 
sur $\C^{*}$ invariantes par $\sq$, ce qui donne la description
\emph{loxodromique} du corps $\M(\Eq) = \M(\C^{*})^{\sq}$ des
fonctions elliptiques: si $q = e^{2 \ii \pi \tau}$, il s'agit
des m\^emes fonctions et des m\^emes corps. Toute fonction
elliptique $f \in \M(\Eq)$ non triviale admet un \emph{diviseur
des z\'eros et des p\^oles sur $\Eq$}, not\'e $\div_{\Eq}(f)$.
En tant que fonction $q$-invariante sur $\C^{*}$, elle admet 
\'egalement un diviseur sur $\C^{*}$, not\'e $\div_{\C^{*}}(f)$. \\

La th\'eorie classique des fibr\'es en droites (ou des diviseurs)
sur $E = \C/(\Z + \Z \tau)$ est la th\'eorie des fonctions Theta
de la forme
$\Theta(\tau,x)$ (\cite{MumfordTheta}): par trivialisation sur
le rev\^etement universel $\C$, on identifie les sections d'un
tel fibr\'e comme des fonctions sur $\C$. La trivialisation 
sur $\C^{*}$ fait de m\^eme apparaitre les \emph{fonctions Theta
de Jacobi}. Nous utiliserons principalement la fonction:
$$
\theta_{q}(z) = \sum_{n \in \Z} q^{-n(n+1)/2} z^{n}.
$$
Cette fonction, qui est holomorphe sur $\C^{*}$ y admet 
la factorisation (\emph{formule du triple produit de Jacobi}):
$$
\theta_{q}(z) = \prod_{n \geq 1} (1 - q^{-n})
\prod_{n \geq 1} (1 + q^{-n} z) \prod_{n \geq 0} (1 + q^{-n} z^{-1}).
$$
Ses z\'eros sont donc les points de $[-1;q]$, compt\'es avec
multiplicit\'e $1$. Comme elle n'a pas de p\^oles, son diviseur
sur $\C^{*}$ est: 
$$
\div_{\C^{*}}(\theta_{q}) = \sum_{a \in [-1;q]} [a].
$$
La fonction $\theta_{q}$ v\'erifie l'\'equation fonctionnelle:
$$
\sq \theta_{q} = z \theta_{q}.
$$
C'est donc une section du fibr\'e en droite $\F_{(z)}$. En tant 
que section, elle admet un diviseur sur $\Eq$:
$$
\div_{\Eq}(\theta_{q}) = [\overline{-1}];
$$
autrement dit, bien que ses valeurs sur $\Eq$ ne soient pas 
d\'efinies, elle y admet le z\'ero simple $\overline{-1}$ et
pas de p\^oles. Nous noterons, pour $a \in \C^{*}$:
$$
\theta_{q,a}(z) = \theta_{q}(z/a).
$$
C'est une fonction holomorphe sur $\C^{*}$, qui y v\'erifie
l'\'equation aux $q$-diff\'erences 
$\sq \theta_{q,a} = \dfrac{z}{a} \theta_{q,a}$
(c'est donc une section de $\F_{(z/a)}$) et l'on a:
$\div_{\Eq}(\theta_{q}) = [\overline{-a}]$. Elle permet
de construire les \emph{$q$-caract\`eres}:
$$
e_{q,a} = \dfrac{\theta_{q}}{\theta_{q,a}} \cdot
$$
On a $\sq e_{q,a} = e_{q,a}$ (c'est donc une section de $\F_{(a)}$)
et $\div_{\Eq}(e_{q,a}) = [\overline{-1}] - [\overline{-a}]$.


\subsubsection{Modules aux $q$-diff\'erences}
\label{subsubsection:modules}

Soit $K$ l'un de nos corps de fonctions, muni de l'automorphisme $\sq$. 
Notre objet est l'\'etude des \emph{\'equations aux $q$-diff\'erences}:
\begin{equation}
\label{equation:eqd}
\sq X = A X, \quad A \in \GL_{n}(K).
\end{equation}
Les cas d'int\'er\^et sont ceux des corps $\Kr$ et $\Ka$. Pour avoir
un bon formalisme alg\'ebrique, on d\'efinit un anneau de polyn\^omes 
de Laurent non commutatifs:
$$
\D = K\left< \sigma,\sigma^{-1}\right>,
$$
par la r\`egle de (non-)commutation: $\sigma z = q z \sigma$. 
Nous noterons $\DM$ la cat\'egorie des $\D$-modules \`a gauche 
de longueur finie. \\

Un objet de $\DM$ peut se r\'ealiser sous la forme $M = (V,\Phi)$ 
o\`u $V$ est un $K$-espace vectoriel de dimension finie et $\Phi$ 
un automorphisme $\sq$-lin\'eaire, \cad\ tel que 
$\Phi(\lambda x) = \sq(\lambda) \Phi(x)$ (l'action de $\sigma$
sur $M$ est alors celle de $\Phi$). Apr\`es choix d'une base,
on peut m\^eme \'ecrire $M = M_{A} = (K^{n},\Phi_{A})$, o\`u
$\Phi_{A}(X) = A^{-1} (\sq X)$ pour une matrice $A \in \GL_{n}(K)$.
Si $A \in \GL_{n}(K)$ et $B \in \GL_{p}(K)$, un morphisme de $A$
dans $B$ est une matrice $F \in \Ma_{p,n}(K)$ telle que:
\begin{equation}
\label{equation:morphisme}
(\sq F) A = B F.
\end{equation}
Si par exemple $F$ est un isomorphisme, alors on a la formule
de transformation de jauge:
$$
B = F[A] = (\sq F) A F^{-1}.
$$
Lorsque par exemple $K = \Kr$, on retrouve l'\'equivalence rationnelle
des \'equations aux $q$-diff\'erences, \'etudi\'ee par Birkhoff. Par
ailleurs, $\DM$ est une cat\'egorie ab\'elienne que l'on peut munir
de constructions tensorielles (par exemple \cite{vdPS} ou \cite{RS2}), 
et il n'est pas tr\`es difficile de v\'erifier que c'est une cat\'egorie
tannakienne (\cite{DF}). En particulier, outre le produit tensoriel,
les constructions suivantes sont disponibles.
\begin{enumerate}
\item{Hom interne: si $M = (V,\Phi)$ et $N = (W,\Psi)$ sont deux modules,
alors $\L_{K}(V,W)$ muni de $f \mapsto \Psi \circ f \circ \Phi^{-1}$
est le module not\'e $\underline{Hom}(M,N)$. On a une adjonction:
$Hom\left(M,\underline{Hom}(M',M'')\right) = Hom(M \otimes M',M'')$.}
\item{Objet unit\'e: c'est le module $\underline{1} = M_{(1)} = (K,\sq)$, 
qui mod\'elise l'\'equation triviale $\sq f = f$. Il est neutre pour 
le produit tensoriel et l'on a $\underline{Hom}(\underline{1},M) = M$.}
\item{Dual: c'est $M^{\vee} = \underline{Hom}(M,\underline{1})$.
Si $M = (K^{n},\Phi_{A})$, on peut le d\'ecrire comme
$M^{\vee} = (K^{n},\Phi_{A^{\vee}})$, o\`u $A^{\vee} = {}^{t} A^{-1}$.}
\item{Foncteur des sections: $\Gamma(M) = Hom(\underline{1},M)$
s'identifie au $K^{\sq}$-espace vectoriel des points fixes de 
$M = (V,\Phi)$ (les $X \in V$ tels que $\Phi(X) = X$). Par exemple
$\Gamma(M_{A})$ est l'espace des solutions de (\ref{equation:eqd})
dans $K$. Le foncteur $\Gamma$ est exact \`a gauche. Son premier
foncteur d\'eriv\'e est $\Gamma^{1}(M) = Ext(\underline{1},M)$.}
\end{enumerate}
Notons, pour un usage futur, les identifications naturelles suivantes:
$\underline{Hom}(M,N) = M^{\vee} \otimes N$ et par cons\'equent:
$Hom(M,N) = \Gamma(M^{\vee} \otimes N)$. Par un argument d'alg\`ebre
homologique, on en d\'eduit $Ext(M,N) = \Gamma^{1}(M^{\vee} \otimes N)$.
Les $Ext^{n}(M,N)$ pour $n \geq 2$ sont nuls, car $\D$ est euclidien 
\`a gauche. \\

Les objets de $\DMr$, $\DMa$ et $\DMf$ sont appel\'es \emph{modules 
aux $q$-diff\'erences}. Dans la pratique, on ne distingue pas toujours
le module $M_{A}$, l'\'equation (\ref{equation:eqd}) et la matrice $A$.
Nous \'etudierons de pr\`es des foncteurs fibres sur ces trois cat\'egories. 
La premi\`ere (``cas global'') est \emph{a priori} notre cat\'egorie 
d'int\'er\^et, mais l'\'etude locale pr\'eliminaire conduit \`a examiner 
$\DMa$ (``cas local analytique'') et $\DMf$ (``cas formel''). \\

Contrairement \`a ce qui se fait pour les \'equations diff\'erentielles,
ni la classification ni la th\'eorie de Galois ne reposent fortement
sur l'\'etude des solutions. La raison est essentiellement celle-ci. 
Pour construire une solution matricielle fondamentale $\X$ de 
l'\'equation (\ref{equation:eqd}), il faut un assez gros corps de
fonctions, mettons $\Kw$. Les solutions vectorielles sont alors les
$X = \X C$, o\`u le vecteur colonne $C$ a ses coefficients 
dans le corps des constantes $\Kw^{\sq} = \M(\Eq)$: c'est un trop
gros corps des constantes (en th\'eorie de Galois diff\'erentielle,
le corps des constantes qui fournit les invariants de classification
est $\C$). Ces raisons et la strat\'egie qui en d\'ecoule ont \'et\'e 
d\'etaill\'ees dans \cite{JSGAL}, \cite{RSZ} et \cite{RS1}. Si l'on
ne tient pas \`a une th\'eorie qui fournisse des invariants transcendants,
alors l'approche alg\'ebrique de van der Put et Singer dans \cite{vdPS}
est appropri\'ee.



\section{Constructions locales}
\label{section:constructionslocales}


\subsection{Construction g\'eom\'etrique g\'en\'erale}
\label{subsection:constructiongeometrique}

Nous noterons d\'esormais $\E = \DMa$ la cat\'egorie des
modules (ou \'equations) aux $q$ diff\'erences sur $\Ka$.
Soit $A(z) \in \GL_{n}(\Ka)$. Soit $D$ un disque \'epoint\'e
en son centre $0$ tel que $A \in \GL_{n}(\O(D))$, \ie\ $A$
et $A^{-1}$ sont holomorphes sur $D$. Sur le fibr\'e trivial
$D \times \C^{n}$ (resp., sur sa base $D$), on d\'efinit une action 
\emph{partielle} de $q^{\Z}$ par l'action de son g\'en\'erateur:
$(z,X) \mapsto (q z,A(z) X)$ (resp. $z \mapsto q z$). (Il reviendrait 
donc au m\^eme de consid\'erer l'action du semi-groupe $q^{-\N}$.) 
Via la projection $D \times \C^{n} \rightarrow D$, ces actions sont
compatibles. On a donc une relation d'\'equivalence $\sim_{A}$ sur 
$D \times \C^{n}$ (resp. $\sim$ sur $D$) engendr\'ee par les relations 
$(z,X) \sim_{A} (q z,A(z) X)$ (resp. $z \sim q z$). On en d\'eduit,
par passage au quotient, un diagramme commutatif:
\begin{equation*}
\begin{CD}
D \times \C^{n} @>pr_{1}>> D \\
@V{}VV                          @VV{\pi}V      \\
\F_{A} = \dfrac{D \times \C^{n}}{\sim_{A}} @>{}>> \Eq = \dfrac{D}{\sim}
\end{CD}
\end{equation*}
Il est en effet bien \'evident que le quotient de la surface 
de Riemann $D$ par $\sim$ est bien la courbe elliptique $\Eq$.
La seule nouveaut\'e ici, par rapport au formalisme g\'en\'eral
de l'introduction, est que la projection $\pi: D \rightarrow \Eq$
n'est plus un rev\^etement, c'est seulement un isomorphisme local.
La ligne du bas d\'ecrit un \emph{fibr\'e vectoriel holomorphe}
$\F_{A}$ sur $\Eq$, et la plus grande partie du discours des
sections \ref{subsubsection:correspondancedeWeil}
et \ref{subsubsection:fibressurEq} se transpose ici. On peut
d\'ecrire le fibr\'e $\F_{A}$ en termes de cocycles, comme dans
\cite{Gunning}: c'est fait dans \cite{AG}. Voici la description
du faisceau des sections. Soit $V$ un ouvert de $\Eq$.
Alors l'espace des sections de $\F_{A}$ sur $V$ est:
\begin{equation}
\label{equation:faisceausections}
\Gamma(V,\F_{A}) = 
\{X \in \O\left(\pi^{-1}(V) \cap D\right)^{n} \mid
\forall z \in \pi^{-1}(V) \cap q^{-1} D \;,\; X(q z) = A(z) X(z) \}.
\end{equation}
(La condition $z \in q^{-1} D$ \'equivaut \`a 
$z \in D \, \wedge \, q z \in D$.) On peut v\'erifier directement
que ce faisceau est localement libre sur $\Eq$ (\cite{Praagman},
\cite{AG}). \\

Si l'on remplace $D$ par $D' \subset D$, on obtient un fibr\'e $\F'_{A}$ 
et, pour tout ouvert $V$ de $\Eq$, un morphisme canonique de restriction
$\Gamma(V,\F_{A}) \rightarrow \Gamma(V,\F'_{A})$ dont il est facile
de v\'erifier qu'il est bijectif. Ainsi, le fibr\'e $\F_{A}$ ne
d\'epend pas du choix particulier du disque $D$. Ce dernier peut
d'ailleurs \^etre remplac\'e par un disque \'epoint\'e topologique.
Dans (\ref{equation:faisceausections}), il faut alors remplacer
la condition $z \in \pi^{-1}(V) \cap q^{-1} D$ par la condition
$z \in \pi^{-1}(V) \cap q^{-1} D \cap D$. Puisque notre construction
ne d\'epend que du \emph{germe} $(D,0)$ de disque \'epoint\'e $D$ 
au voisinage de $0$, nous noterons plus intrins\`equement:
$$
\F_{A} = \dfrac{(D,0) \times \C^{n}}{\sim_{A}} \longrightarrow 
\Eq = \dfrac{(D,0)}{\sim} \cdot
$$

On obtient ainsi un foncteur $A \leadsto \F_{A}$ de $\E = \DMa$ dans 
la cat\'egorie des fibr\'es vectoriels holomorphes sur $\Eq$. 
Les propri\'et\'es (faciles) suivantes sont alors essentielles 
pour la th\'eorie de Galois: \emph{ce foncteur est exact, fid\`ele
et $\otimes$-compatible}. Dans la terminologie de \cite{DF}, on
dit que c'est un foncteur fibre de $\E$ sur $\Eq$. Il permet 
de construire des familles de foncteurs fibres sur $\C$, respectivement
index\'ees par $\Eq$ et par $\C^{*}$: voir \cite{RS1} et \cite{RS2}. \\

D'un point de vue purement fonctoriel, nous verrons que le fonteur
$A \leadsto \F_{A}$ est essentiellement surjectif; mais cela ne sert
\`a rien, car il n'est pas pleinement fid\`ele. Par exemple, il peut
associer des fibr\'es isomorphes \`a des modules qui ne le sont pas
(voir un exemple section \ref{subsubsection:HarderNarasimhan}). 
Cette question sera abord\'ee \`a la la section 
\ref{subsubsection:fibrespentesentieres}. \\

Pratiquement, le probl\`eme se pr\'esente ainsi. Soient 
$A \in \GL_{n}(\Ka)$ et $B \in \GL_{p}(\Ka)$ les matrices de deux objets 
$M_{A},M_{B}$ de $\E$. Soit $\phi: \F_{A} \rightarrow \F_{B}$ un
morphisme. D'apr\`es les g\'en\'eralit\'es de l'introduction, on
peut d\'ecrire $\phi$ comme une application holomorphe $F$ de $D$ dans
$\Ma_{p,n}(\C)$ qui satisfait l'\'equation (\ref{equation:morphisme}).
Pour en faire un morphisme dans $\E$, il faudrait le prolonger
m\'eromorphiquement en $0$. \emph{Mais on ne sait (en g\'en\'eral) 
rien du mode de croissance de $F$ en $0$}. \\

\Exs
1. Un morphisme de $\F_{(1)}$ dans $\F_{(z)}$ s'identifie \`a une fonction
holomorphe $f: \C^{*} \rightarrow \C$ telle que $f(q z) = z f(z)$,
autrement dit, \`a une section holomorphe de $\F_{(z)}$. (Ainsi, 
$\F_{(1)}$ se comporte comme l'objet unit\'e $\underline{1}$ d\'ecrit
en \ref{subsubsection:modules}.) Par exemple, $\theta_{q}$ r\'ealise
un tel morphisme, et il a une singularit\'e essentielle en $0$.
D'ailleurs, il n'existe aucun morphisme non trivial de 
$\underline{1} = M_{(1)}$ dans $M_{(z)}$: en effet, ce serait
une s\'erie de Laurent $f \in \Ka$ telle que $\sq f = z f$,
ce qui implique $f = 0$. (Cela reste vrai dans $\DMf$.) \\
2. Soient $A = \begin{pmatrix} 1 & 0 \\ 0 & z \end{pmatrix}$ et
$F = \begin{pmatrix} 1 & 0 \\ \theta_{q} & 1 \end{pmatrix}$. Alors
$F$ r\'ealise un automorphisme de $\F_{A}$, qui ne provient pas
d'un automorphisme de $M_{A}$. \\
3. Soit $A = \begin{pmatrix} 1 & 0 \\ 0 & z^{2} \end{pmatrix}$.
Nous verrons en \ref{subsubsection:HarderNarasimhan} que $\F_{A}$
est somme de deux fibr\'es en droites de degr\'e $1$, alors que
$M_{A}$ est ind\'ecomposable. \\

\Rems
1. Le fibr\'e $\F_{A}$ d\'efini \`a partir du quotient de $D \times \C^{n}$ 
se rel\`eve par $\pi: \C^{*} \rightarrow \Eq$ en un fibr\'e trivial sur 
$\C^{*}$ (section \ref{subsubsection:fibressurEq}). Il est donc de la forme
$\F_{A'}$, o\`u $A'$ est holomorphe de $\C^{*}$ dans $\GL_{n}(\C)$.
Le fait qu'il s'agisse de deux rel\`evements du m\^eme fibr\'e signifie
qu'il existe $F: D \rightarrow \GL_{n}(\C)$ holomorphe tel que
$(\sq F) A = A' F$. Les matrices $A$ et $A'$ sont holomorphiquement
\'equivalentes sur $D$, mais $A$ se prolonge m\'eromorphiquement en $0$
alors que $A'$ se prolonge holomorphiquement \`a $\C^{*}$. \\
2. Tout fibr\'e sur une surface de Riemann compacte est m\'eromorphiquement
trivial (\cite{Gunning}, p. 103)
\footnote{Gunning l'affirme en g\'en\'eral, en faisant r\'ef\'erence
\`a la page 43 o\`u c'est prouv\'e seulement pour la droite projective; 
en fait, la d\'emonstration s'adapte sans probl\`eme.}. 
Il existe donc $F: D \rightarrow \GL_{n}(\C)$ m\'eromorphe tel que
$A = F[I_{n}]$: $F$ est donc une solution matricielle fondamentale 
\emph{m\'eromorphe } de (\ref{equation:eqd}). (C'est en substance 
l'argument de Praagman dans \cite{Praagman}.)


\subsubsection{\'Equations fuchsiennes}
\label{subsubsection:equationsfuchsiennes}

Nous dirons que le module $M$ de $\E$ ou $\DMf$ est \emph{fuchsien}
s'il est de la forme $M_{A}$, o\`u $A(0) \in \GL_{n}(\C)$. De m\^eme,
si $K = \Ka$ ou $\Kf$, l'\'equation (\ref{equation:eqd}) est dite
\emph{fuchsienne} si elle est \'equivalente \`a une \'equation de
matrice $B$ telle que $B(0) \in \GL_{n}(\C)$. Dans le cas d'un module 
de $\DMr$ ou d'une \'equation \`a coefficients dans $\Kr$, consid\'er\'es
via le plongement $\Kr \hookrightarrow \Ka$, on dit \emph{fuchsien(ne)
en $0$}. Cette propri\'et\'e \'equivaut \`a des conditions de croissance
des solutions (\cite{JSAIF}) ou de polygone de Newton (section 
\ref{subsubsection:polNewt}). Pour une caract\'erisation plus intrins\`eque, 
en termes de ``r\'eseau stable'', voir \cite{vdPS}, \cite{LDV}. \\

Un lemme cl\'e dit alors que \emph{toute \'equation fuchsienne est
localement \'equivalente \`a une \'equation \`a coefficients constants},
autrement dit, il existe $A^{(0)} \in \GL_{n}(\C)$ et $F \in \GL_{n}(\Ka)$
tels que $A = F\left[A^{(0)}\right]$. Comme dans le cas des \'equations
diff\'erentielles, se lemme se prouve en deux \'etapes: \'elimination
des \emph{r\'esonnances} \`a l'aide de transformations de shearing;
d\'etermination d'un unique $F$ formel tangent \`a l'identit\'e, et
preuve de convergence. \\

La cat\'egorie $\Ef$ des \'equations fuchsiennes en $0$ est,
par d\'efinition, la sous-cat\'egorie pleine de $\E$ form\'ee
des objets fuchsiens. (En fait, on constate \emph{a posteriori}
que l'on peut aussi bien partir de $\DMr$ ou de $\DMf$.) C'est
une sous-cat\'egorie tannakienne de $\E$ (elle est stable par
toutes les op\'erations lin\'eaires). La sous-cat\'egorie pleine
form\'ee des \'equations \`a coefficients constants est \'egalement
une sous-cat\'egorie tannakienne, que nous noterons $\P$. Il d\'ecoule
du lemme-cl\'e que \emph{l'inclusion $\P \hookrightarrow \Ef$ est
une \'equivalence de cat\'egories.} \\

Pour tout objet $M_{A}$ de $\P$ d\'efini par $A \in \GL_{n}(\C)$, 
le fibr\'e $\F_{A}$ admet la construction simplifi\'ee: 
$$
\F_{A} = \dfrac{\C^{*} \times \C^{n}}{(z,X) \sim (qz,AX)} \cdot
$$
C'est donc un fibr\'e plat.

\begin{thm}[\cite{JSGAL}]
\label{thm:pleinefidelite}
Le foncteur $A \leadsto \F_{A}$ de $\P$ dans la cat\'egorie 
$Fib_{p}(\Eq)$ des fibr\'es plats sur $\Eq$ est une \'equivalence 
de cat\'egories tannakiennes.
\end{thm}
\Pr
Soit $\phi: \F_{A} \rightarrow \F_{B}$, repr\'esent\'e par une matrice
$F$ holomorphe sur $\C^{*}$ et telle que $(\sq F) A = B F$. On a donc 
$F(q z) = B F(z) A^{-1}$, d'o\`u l'on d\'eduit facilement que $F$ a une 
croissance mod\'er\'ee en $0$, donc un prolongement m\'eromorphe, d'o\`u 
la pleine fid\'elit\'e. (En fait, il n'est pas tr\`es difficile de voir
que $F$ est \`a coefficients dans $\C[z,z^{-1}]$, \cf\ \cite{JSGAL}.) \\
Pour l'essentielle surjectivit\'e, on part d'un fibr\'e plat d\'efini par
la repr\'esentation telle que $1 \mapsto A_{1}$, $\tau \mapsto A_{\tau}$.
Les matrices $A_{1}$ et $A_{\tau}$ commutent. Il existe donc un logarithme 
$2 \ii \pi B$ de $A_{1}$ qui commute avec $A_{1}$ et $A_{\tau}$. En posant
$G(x) = e^{2 \ii \pi x B}$, on voit que $G(x+1) = A_{1} G(x)$ et
$G(x+\tau) A = A_{\tau} G(x)$, o\`u $A = A_{\tau} A_{1}^{-\tau}$.
Le fibr\'e est donc isomorphe \`a $\F_{A}$.
\hfill $\Box$ \\

Au fibr\'e plat $\F_{A}$ est associ\'ee par la correspondance de Weil
la repr\'esentation de $\pi_{1}(\Eq) = \Z + \Z \tau$ d\'efinie par 
$1 \mapsto I_{n}$, $\tau \mapsto A$. Le th\'eor\`eme dit que toute
repr\'esentation de $\pi_{1}(\Eq)$ est \'equivalente (dans ce sens \'etendu) 
\`a une repr\'esentation de cette forme. La cat\'egorie $Fib_{p}(\Eq)$ 
(dont les morphismes sont \emph{tous} les morphismes de fibr\'es vectoriels 
holomorphes) est donc \'equivalente \`a la cat\'egorie dont les objets sont 
les repr\'esentations du groupe $\Aut(\C^{*} / \Eq) = q^{\Z} \simeq \Z$, et 
les morphismes sont les morphismes \'equivariants de la section 
\ref{subsubsection:correspondancedeWeil}. La cat\'egorie $\Ef$ 
est donc \'equivalente \`a la cat\'egorie $Rep_{\Eq}(\Z)$ obtenue 
en \'epaississant la cat\'egorie $Rep(\Z)$ des repr\'esentations
(complexes de dimension finie) de $\Z$. Dans \cite{JSGAL} (et,
par une autre voie, dans \cite{BG}) on en d\'eduit la description
du groupe de Galois de $\Ef$. L'action de ce groupe de Galois est
explicit\'ee dans \cite{RS1} et \cite{RS2}.


\subsubsection{Quelques exemples de rang $1$}
\label{subsubsection:rang1}

Un module de rang $1$ dans $\E$ est de la forme $M_{(a)}$,
o\`u $a \in \Ka^{*}$. \'Ecrivons $a = a_{0} z^{\mu} u$, o\`u
$a_{0} \in \C^{*}$, $\mu = v_{0}(a) \in \Z$ (valuation $z$-adique 
de $a$) et $u = 1 + u_{1} z + \cdots$. Pour construire des
solutions, on applique la section \ref{subsubsection:fonctions}. \\

L'\'equation $\sq f = u f$ admet la solution r\'eguli\`ere
$v(z) = \prod\limits_{n \geq 1} u(q^{-n} z)$. Il revient au
m\^eme de dire que $v: (1) \rightarrow (u)$ est une \'equivalence
analytique, ou que $v: \underline{1} \rightarrow M_{(u)}$ est
un isomorphisme. Le fibr\'e $\F_{(u)}$ est donc isomorphe au 
fibr\'e en droites trivial $\underline{1} = \F_{(1)}$. \\

L'\'equation $\sq f = z^{\mu} f$ admet la solution $\theta_{q}^{\mu}$.
Le module $M_{(z^{\mu})}$ est isomorphe \`a $M_{(z)}^{\otimes \mu}$,
puissance tensorielle $\mu$\ieme de $M_{(z)}$. (Comme tout module 
de rang $1$, il admet son dual comme inverse pour $\otimes$, ses 
puissances tensorielles n\'egatives sont donc d\'efinies.) 
Le fibr\'e associ\'e est $\F_{(z^{\mu})} \simeq \F_{(z)}^{\otimes \mu}$,
dont une section est $\theta_{q}^{\mu}$, de diviseur
$\div_{\Eq} \theta_{q}^{\mu} = \mu [\overline{-1}]$. \\

L'\'equation fuchsienne scalaire $\sq f = a_{0} f$ admet pour solution 
le $q$-caract\`ere $e_{q,a_{0}}$, qui est une section du fibr\'e plat
$\F_{a_{0}}$. Comme 
$\div_{\Eq}(e_{q,a_{0}}) = [\overline{-1}] - [\overline{-a_{0}}]$,
le degr\'e de ce fibr\'e est $0$. (Le lecteur remarquera que, si
$a_{0} \in q^{\Z}$, le degr\'e est nul et le module et le fibr\'e 
sont en fait triviaux.) \\

Le module $M_{(a)}$ (resp. le fibr\'e $\F_{(a)}$ associ\'e) est isomorphe 
au produit tensoriel de ces trois modules (resp. de ces trois fibr\'es).
Le degr\'e de $\F_{(a)}$ est donc $\mu$. Sa pente (quotient du degr\'e par 
le rang, \cf\ \cite{Seshadri}) est donc $\mu$.


\subsection{\'Equations irr\'eguli\`eres}
\label{subsection:equationsirregulieres}


\subsubsection{Polygone de Newton}
\label{subsubsection:polNewt}

On montre que tout module aux $q$-diff\'erences sur $K = \Kr$, $\Ka$ 
ou $\Kf$ peut se mettre sous la forme  $M = \D/\D P$ (``lemme du vecteur 
cyclique'' et euclidianit\'e \`a gauche de $\D$). On peut prendre $P$
sous la forme $a_{n} \sigma^{n} + \cdots + a_{0}$, o\`u $a_{0} a_{n} \neq 0$.
Notons $v_{0}(a)$ la valuation $z$-adique de $a \in K$. La fronti\`ere de 
l'enveloppe convexe de $\{(i,j) \in \N \times \Z \mid 
0 \leq i \leq n \text{~et~} j \geq v_{0}(a_{j})\}$
est form\'ee de deux demi-droites verticales et de $k$ vecteurs
$(r_{1},d_{1}),\ldots,(r_{k},d_{k}) \in \N^{*} \times \Z$. Notant
$\mu_{i} = \dfrac{d_{i}}{r_{i}} \in \Q$, on indexe ces vecteurs de gauche 
\`a droite, de sorte que $\mu_{1} < \cdots < \mu_{k}$. On prouve que 
les \emph{pentes} 
\footnote{Depuis \cite{RS1}, \cite{RS2}, nous avons adopt\'e pour les
pentes une convention \emph{oppos\'ee} \`a celle qui pr\'evalait dans 
\cite{JSFIL}, \cite{RSZ}, \cite{JSSTO}.}
$\mu_{i}$ et leurs \emph{multiplicit\'es} $r_{i}$ 
ne d\'ependent que de $M$. \\

\Exs
1. Soit $L = \sigma - a$ o\`u $a \in \Ka^{*}$. L'\'equation $L f = 0$, 
\cad\ $\sq f = a f$, est mod\'elis\'ee par le module $M_{(a)}$, dont 
on peut d\'emontrer qu'il est isomorphe \`a $\D / \D \L^{\vee}$, o\`u 
$L^{\vee} = a \sigma - 1$. Les pentes de $M_{(a)}$ se calculent
avec $L^{\vee}$: l'unique pente est $\mu = v_{0}(a)$. Le fibr\'e $\F_{(a)}$ 
est de rang $1$ et de degr\'e $\mu$ (section \ref{subsubsection:rang1}),
donc de pente $\mu$. \\
\label{Ex:Tschakaloff}
2. Soit $L = q z \sigma^{2} - (1+z) \sigma + 1 = (\sigma - 1)(z \sigma - 1)$.
L'\'equation $L f$ est int\'eressante, parce qu'elle est satisfaite par
la \emph{s\'erie de Tschakaloff} $\sum\limits_{n \geq 0}  q^{n(n-1)/2} z^{n}$,
qui est un $q$-analogue de la s\'erie d'Euler. On en fait une \'equation
vectorielle de type (\ref{equation:eqd}) en posant (par exemple)
$X = \begin{pmatrix} f \\ z \sq f - f \end{pmatrix}$ et
$A = \begin{pmatrix} z^{-1} & z^{-1} \\ 0 & 1 \end{pmatrix}$.
Le module $M_{A}$ est isomorphe \`a \`a $\D / \D \L^{\vee}$, o\`u 
$L^{\vee} = (\sigma - z)(\sigma - 1) = \sigma^{2} - (1+z) \sigma + z$.
Ses pentes sont $-1$ et $0$. La forme triangulaire de la matrice
indique qu'il y a une suite exacte:
$0 \rightarrow M_{(z^{-1})} \longrightarrow M_{A} \longrightarrow
M_{(1)} \rightarrow 0$,
o\`u l'inclusion a pour matrice $\begin{pmatrix} 1 \\ 0 \end{pmatrix}$
et la projection $\begin{pmatrix} 0 & 1 \end{pmatrix}$. Il y a donc
\'egalement une suite exacte de fibr\'es:
$0 \rightarrow \F_{(z^{-1})} \longrightarrow \F_{A} \longrightarrow
\F_{(1)} \rightarrow 0$. \\

La donn\'ee du \emph{polygone de Newton} $N_{M}$ de $M$ \'equivaut
\`a celle de la fonction $r_{M}$ de $\Q$ dans $\N$ telle que 
$\mu_{i} \mapsto r_{i}$ et qui est nulle par ailleurs. Selon \cite{JSFIL}, 
le polygone de Newton est additif pour les suites exactes, multiplicatif 
pour le produit tensoriel et tel que $r_{M^{\vee}}(\mu) = r_{M}(-\mu)$. 
Ces r\`egles sont assez diff\'erentes de celles qui r\'egissent les 
\'equations diff\'erentielles. \\

Nous dirons qu'un module est \emph{pur isocline} s'il admet une seule
pente et \emph{pur} s'il est somme directe de modules purs isoclines
\footnote{Nous avons adopt\'e cette terminologie depuis \cite{RS1}, 
\cite{RS2} (anciens termes: ``pur'', et ``mod\'er\'ement irr\'egulier'';
nos modules purs sont les ``split modules'' de \cite{vdPR}).}.
Les modules fuchsiens sont les modules purs isoclines de pente $0$.


\subsubsection{Filtration par les pentes}
\label{subsubsection:filtration}

\begin{thm}[\cite{JSFIL}]
\label{thm:filtration}
(i) Tout module de $\E$ ou $\DMf$ admet une unique filtration
croissante $(M_{\leq \mu})_{\mu \in \Q}$ telle que les quotients
$M_{(\mu)} = M_{\leq \mu}/M_{< \mu}$ sont purs de pente $\mu$.
(Le rang de $M_{(\mu)}$ est donc $r_{M}(\mu)$.) \\
(ii) la filtration est strictement fonctorielle: le foncteur gradu\'e 
associ\'e: $M \leadsto \gr M = \bigoplus\limits_{\mu \in \Q} M_{(\mu)}$
est exact. Il est de plus fid\`ele et $\otimes$-compatible. \\
(iii) Dans $\DMf$, le foncteur $\gr$ est isomorphe au foncteur identit\'e.
\end{thm}

En termes d'op\'erateurs, ce th\'eor\`eme dit que, pour $K = \Ka$, 
tout $L \in \D$ admet une factorisation $L =L_{1} \cdots L_{k}$ 
o\`u les $L_{i}$ sont purs de pentes $\mu_{1} < \cdots < \mu_{k}$.
Pour $K = \Kf$, une telle factorisation est possible avec un ordre
des pentes arbitraires. En termes de matrices, toute $A \in \GL_{n}(\Ka)$
peut se mettre sous la forme triangulaire sup\'erieure par blocs:
$$
A = \begin{pmatrix}
A_{1}  & \ldots & \ldots & \ldots & \ldots \\
\ldots & \ldots & \ldots  & U_{i,j} & \ldots \\
0      & \ldots & \ldots   & \ldots & \ldots \\
\ldots & 0 & \ldots  & \ldots & \ldots \\
0      & \ldots & 0       & \ldots & A_{k}    
\end{pmatrix},
$$
o\`u chaque $A_{i}$ est pure de pente $\mu_{i}$. Si $A \in \GL_{n}(\Kf)$,
on peut prendre les $U_{i,j} = 0$. \\

\Ex
Soit $L = q z \sigma^{2} - (1+z) \sigma + 1$, dont les pentes sont 
$0$ et $-1$. L'existence analytique de la filtration vient de la
factorisation analytique $L  = (\sigma - 1)(z \sigma - 1)$ (qui
remonte \`a Adams, Birkhoff et Guenther). L'op\'erateur dual
se factorise: $L^{\vee} = (\sigma - z)(\sigma - 1)$, d'o\`u
une suite exacte pour le module associ\'e $M = \D / \D \L^{\vee}$:
$$
0 \rightarrow \D/\D(\sigma - z) \longrightarrow M
\longrightarrow \D/\D(\sigma - 1) \rightarrow 0.
$$
C'est cette suite exacte qui permet de construire la forme 
triangulaire $A = \begin{pmatrix} z^{-1} & z^{-1} \\ 0 & 1 \end{pmatrix}$
obtenue page \pageref{Ex:Tschakaloff}. En effet, le module
$\D/\D(\sigma - z)$ se d\'ecrit \'egalement comme $(\Ka, \Phi_{z^{-1}})$,
et correspond donc \`a la matrice $(z^{-1}) \in \GL_{1}(\Ka)$.
De m\^eme, $\D/\D(\sigma - 1)$ correspond \`a $(1) \in \GL_{1}(\Ka)$. 

\begin{cor}
Le fibr\'e $\F_{A}$ est muni d'une filtration 
$\F_{1} \subset \cdots \subset \F_{k}$ telle que chaque $\F_{i}/\F_{i-1}$
est le fibr\'e associ\'e \`a une \'equation pure de pente $\mu_{i}$.
\end{cor}

\Rem
En fait, sous cette forme, l'\'enonc\'e ci-dessus est presque vide.
Selon le th\'eor\`eme 10, p. 63 de \cite{Gunning}, tout fibr\'e 
sur $\Eq$ peut \^etre d\'ecrit par une matrice triangulaire 
sup\'erieure (sans blocs). Avec les calculs de la section
\ref{subsubsection:rang1}, on peut obtenir des termes diagonaux
de la forme $\alpha_{i} z^{d_{i}}$ ($\alpha_{i} \in \C^{*}$,
$d_{i} \in \Z$). De plus, chaque fois que $i < j$ et $d_{i} > d_{j}$,
on peut permuter les deux termes diagonaux correspondants. 
En effet, cela se ram\`ene \`a la constatation que l'\'equation
$\alpha_{j} z^{d_{j}} \sq f - \alpha_{i} z^{d_{i}} f = u$
admet, pour tout $u \in \O(\C^{*})$, une solution $f \in \O(\C^{*})$;
et ce point est imm\'ediat par identification des s\'eries de Laurent.
On obtient donc un r\'esultat apparemment analogue au th\'eor\`eme
mais plus fort pour les fibr\'es. Cependant, les exposants $d_{i}$
qui apparaissent ici n'ont aucune signification intrins\`eque en
termes de fibr\'es. (Pour une mise en forme plus intris\`eque,
voir la section \ref{subsubsection:HarderNarasimhan}.) \\

\`A partir de ce th\'eor\`eme, deux voies distinctes ont \'et\'e suivies.
En supposant les pentes enti\`eres, on a une description simple des blocs 
purs $A_{i}$ (section \ref{subsection:pentesentieres}). 
On en tire, \emph{par voie transcendante}, des cons\'equences 
pour la classification (\cite{RSZ}, \cite{JSSTO}, section 
\ref{subsection:classification}) et pour la th\'eorie de Galois 
(\cite{RS1}, \cite{RS2}, section \ref{subsection:calculsdinvariants}). \\

R\'ecemment, van der Put et Reversat ont \'elucid\'e la structure
des modules purs dans le cas g\'en\'eral (\cite{vdPR}). Dans le cas
d'un module ind\'ecomposable de rang $r$ et pente $\mu = \dfrac{d}{r}$,
on a $d \wedge r = 1$ et la description est similaire \`a celle des
fibr\'es ind\'ecomposables sur $\Eq$ par Atiyah (\cite{Atiyah}, 
\cite{Koszul}), qu'elle permet de retrouver de mani\`ere plus simple.
L'extension des r\'esultats de Ramis, Sauloy et Zhang au cas des pentes
rationnelles \`a l'aide de \cite{vdPR} ne semble pas avoir \'et\'e
faite. Elle devrait entrainer des complications alg\'ebriques, mais
peut-\^etre pas mettre en jeu d'id\'ees analytiques nouvelles.
Dans \emph{loc. cit.} , van der Put et Reversat enrichissent de plus
le fibr\'e d'une connexion m\'eromorphe, mais cela ne modifie pas 
le probl\`eme de la pleine fid\'elit\'e. Nous proposerons une autre 
structure \`a la section \ref{subsubsection:fibrespentesentieres}.


\subsection{\'Equations irr\'eguli\`eres \`a pentes enti\`eres}
\label{subsection:pentesentieres}


\subsubsection{Description des \'equations \`a pentes enti\`eres}
\label{subsubsection:equationspentesentieres}

Pour tout $d \in \N^{*}$, la sous-cat\'egorie pleine de $\E = \DMa$
form\'ee des \'equations \`a pentes dans $\dfrac{1}{d} \Z$ est une
sous-cat\'egorie tannakienne. Pour $d = 1$, on obtient la cat\'egorie
$\Ee$ des \'equations \`a pentes enti\`eres. La sous-cat\'egorie pleine
de $\Ee$ form\'ee des \'equations pures \`a pentes enti\`eres en est 
encore une sous-cat\'egorie tannakienne $\Epe$. \\

Si $M$ est pur (isocline) de pente $\mu \in \Z$, alors 
$M_{(z^{-\mu})} \otimes M$ est fuchsien, donc de la forme $M_{A}$
avec $A \in \GL_{n}(\C)$. On a donc $M \simeq M_{z^{\mu} A}$.
On peut alors am\'eliorer la forme triangulaire des matrices
en une \emph{forme standard}:
\begin{equation} 
\label{equation:formestandard}
A = 
\begin{pmatrix}
z^{\mu_{1}} A_{1}  & \ldots & \ldots & \ldots & \ldots \\
\ldots & \ldots & \ldots  & U_{i,j} & \ldots \\
0      & \ldots & \ldots   & \ldots & \ldots \\
\ldots & 0 & \ldots  & \ldots & \ldots \\
0      & \ldots & 0       & \ldots & z^{\mu_{k}} A_{k}    
\end{pmatrix},
\end{equation}
o\`u les pentes $\mu_{1} < \cdots < \mu_{k}$ sont enti\`eres,
$r_{i} \in \N^{*}$, $A_{i} \in \GL_{r_{i}}(\C)$ ($i = 1,\ldots,k$) 
(ces $\mu_{i}$ et $r_{i}$ constituent le polygone de Newton de $A$) et:
$$
U = (U_{i,j})_{1 \leq i < j \leq k} \in 
\underset{1 \leq i < j \leq k}{\prod} \Ma_{r_{i},r_{j}}(\Ka).
$$
On peut m\^eme imposer aux $U_{i,j}$ d'\^etre polynomiaux
(\cite{RSZ}). \\

Soit $B$ une matrice \`a pentes enti\`eres d\'ecrite de mani\`ere
similaire: blocs diagonaux $z^{\nu_{i'}} B_{i'}$, o\`u
$B_{i'} \in \GL_{s_{i'}}(\C)$ ($i' = 1,\ldots,k'$)
et surdiagonaux $V_{i',j'} \in \Ma_{s_{i'},s_{j'}}(\Ka)$. 
La filtration \'etant fonctorielle, tout morphisme
$F: A \rightarrow B$ est triangulaire sup\'erieur par blocs.
Plus pr\'ecis\'ement, $F$ admet une d\'ecomposition en blocs
$F_{i',i} \in \Ma_{s_{i'},r_{i}}(\Ka)$, nuls pour $i < i'$ et 
tels que, pour $i' \geq i$:
$$
(\sq F_{i',i}) z^{\mu_{i}} A_{i} + 
\sum_{i' \leq l < i} (\sq F_{i',l}) U_{l,i} =
z^{\nu_{i'}} B_{i'} F_{i',i} + \sum_{i' \leq l < i} V_{i',l} F_{l,i}.
$$


\subsubsection{Fibr\'es associ\'es aux \'equations \`a pentes enti\`eres}
\label{subsubsection:fibrespentesentieres}

Pour tout module pur isocline $M \simeq M_{z^{\mu} A}$,
le fibr\'e $\F_{M}$ est isomorphe \`a $\F_{(z^{\mu})} \otimes \F_{A}$.
Disons qu'un fibr\'e est \emph{pur isocline de pente $\mu$} si c'est
le produit tensoriel d'un fibr\'e en droites de degr\'e $\mu$ par
un fibr\'e plat. Il r\'esulte du th\'eor\`eme \ref{thm:filtration}
et de son corollaire que l'on peut associer \`a tout objet $M$ de $\Ee$
un fibr\'e $\F_{M}$ muni d'une filtration \`a quotients purs isoclines 
$\F_{1} \subset \cdots \subset \F_{k}$.

\begin{thm}
Le foncteur $M \leadsto (F_{M},\F_{1} \subset \cdots \subset \F_{k})$
est exact, pleinement fid\`ele et $\otimes$-compatible.
\end{thm}
\Pr
Les objets de la cat\'egorie d'arriv\'ee sont les couples 
$\left(\F,(\F_{\leq \mu})_{\mu \in \Z}\right)$ form\'es d'un fibr\'e 
et d'une filtration croissante par des sous-fibr\'es telle que les 
$\F_{(\mu)} = \F_{\leq \mu}/\F_{< \mu}$ sont des fibr\'es purs 
isoclines de pente $\mu$. Les morphismes de cette cat\'egorie
sont les morphismes de fibr\'es $\phi: \F \rightarrow \F'$ qui 
respectent la filtration: $\phi(\F_{\leq \mu}) \subset \F'_{\leq \mu}$.
La structure tensorielle est d\'efinie en munissant $\F \otimes \F'$
de la filtration:
$$
(\F \otimes \F')_{\leq \mu} = \Im
\left(\sum_{\nu + \nu' \leq \mu} \F_{\leq \nu} \otimes \F'_{\leq \nu'}
\rightarrow \F \otimes \F'\right).
$$
Le seul point nouveau dans l'\'enonc\'e du th\'eor\`eme est que 
ce foncteur est \emph{pleinement} fid\`ele. Soient $A$, $B$ des
objets de $\Ee$ et $\phi: \F_{A} \rightarrow \F_{B}$ un morphisme.
On \'ecrit $\phi$ sous la forme $F \in \Ma_{p,n}\left(\O(\C^{*})\right)$
tel que $(\sq F) A = B F$. L'hypoth\`ese sur le respect des filtrations
dit que $F$ est triangulaire par blocs au sens de la section
pr\'ec\'edente, et que ces blocs v\'erifient les m\^emes \'equations
que ceux des morphismes dans $\Ee$. Le lemme ci-dessous entraine
alors, par r\'ecurrence, que $F$ admet un prolongement m\'eromorphe 
en $0$. (La r\'ecurrence est amorc\'ee par les blocs diagonaux de $F$
qui se d\'eduisent du cas fuchsien.)
\hfill $\Box$ 

\begin{lemma}
Soient $\mu < \nu$, $C \in \GL_{r}(\C)$, $D \in \GL_{s}(\C)$
et $U \in \Ma_{s,r}(\Ka)$. Soit $F \in \Ma_{s,r}\left(\O(\C^{*})\right)$ 
tel que $(\sq F) (z^{\mu} C) - (z^{\nu} D) F = U$.
Alors $F$ est m\'eromorphe en $0$.
\end{lemma}
\Pr
On \'ecrit le d\'eveloppement en s\'erie de Laurent:
$F = \sum\limits_{n \in \Z} F_{n} z^{n}$. On a donc:
$q^{n-\mu} F_{n - \mu} C - D F_{n-\nu} = U_{n}$ pour tout $n \in \Z$.
Puisque $U$ est m\'eromorphe en $0$, il existe donc $n_{0}$ tel que, 
pour $n \geq n_{0}$, on ait: $q^{n} F_{n} = D F_{n - \delta} C^{-1}$,
o\`u $\delta = \nu - \mu \geq 1$. Si l'on n'a pas $F_{n} = 0$ pour
$n \ll 0$, alors $F_{n}$ est de l'ordre de grandeur de $q^{n^{2}/2 \delta}$
lorsque $n \to - \infty$, contredisant la convergence de la s\'erie
de Laurent.
\hfill $\Box$


\subsubsection{Lien avec la filtration de Harder-Narasimhan}
\label{subsubsection:HarderNarasimhan}

La filtration introduite \`a la section pr\'ec\'esente pr\'esente
des ressemblances formelles avec la filtration de Harder-Narasimhan 
(\cite{Seshadri}), mais nous allons voir qu'elles ne sont pas li\'ees.
Cependant, il y a d'intrigantes questions de stabilit\'e. Les calculs
qui suivent sont en grande partie tir\'es de \cite{AG}.

\paragraph{Une famille g\'en\'eriquement stable.}

Soit $A_{u} = \begin{pmatrix} 1 & u \\ 0 & z \end{pmatrix}$, avec
$u \in \Ka$. D'apr\`es \cite{RSZ}, on peut ramener $A_{u}$ \`a la
forme $A_{v}$ avec $v \in \C$ via un unique morphisme de la forme
$\begin{pmatrix} 1 & f \\ 0 & 1 \end{pmatrix}$, o\`u $f \in \Ka$, 
\ie\ trivial sur les gradu\'es associ\'es. \\
Soient maintenant $u,v \in \C$. D'apr\`es la fonctorialit\'e de
la filtration et ce que l'on sait sur les morphismes entre objets
fuchsiens, tout morphisme de $A_{u}$ dans $A_{v}$ est de la forme
$\begin{pmatrix} 1 & f \\ 0 & 1 \end{pmatrix}
\begin{pmatrix} \alpha & 0 \\ 0 & \beta \end{pmatrix}$, o\`u
$f \in \Ka$ et $\alpha,\beta \in \C^{*}$. Or, l'\'equation
correspondante $z \sq f - f = v - \dfrac{\alpha}{\beta} u$
admet pour seule solution formelle 
$f = \left(v - \dfrac{\alpha}{\beta} u\right) \hat{\phi}$, 
o\`u $\hat{\phi}$ est 
la s\'erie de Tschakaloff, et $f$ n'est donc convergente que si
$v = \dfrac{\alpha}{\beta} u \cdot$ Ainsi, $A_{u}$ n'est isomorphe
\`a $A_{v}$ que si $u = v = 0$ ou bien si $u \neq 0, v \neq 0$. 
Les $A_{u}$, 
autrement dit, les extensions de $(z)$ par $(1)$ se r\'epartissent
en deux classes
\footnote{Il s'agit des classes d'isomorphie des modules qui sont 
des extensions de $(z)$ par $(1)$, et non des classes d'extension au 
sens habituel. Ces derni\`eres forment l'espace $Ext\left((z),(1)\right)$,
qui est de dimension $1$ (\cite{RSZ}) et param\'etr\'e par $u$.}: 
scind\'ee ou non. Dans le cas non scind\'e, $A_{u}$ est de plus
ind\'ecomposable (cela se d\'eduit de la fonctorialit\'e de la
filtration par les pentes). \\

Discutons maintenant la \emph{stabilit\'e} du fibr\'e $\F_{A_{u}}$
(\cite{Gunning}, \cite{Seshadri}). Le degr\'e de $\F_{A_{u}}$ est $1$
(extension de $\F_{(z)}$ par $\F_{(1)}$), son rang est $2$, sa pente
est donc $\dfrac{1}{2} \cdot$ Pour qu'il soit instable, il faut, et
il suffit, qu'il admette un sous-fibr\'e en droite de pente, donc de
degr\'e $> \dfrac{1}{2}$, donc $\geq 1$, autrement dit, qu'il admette
une section m\'eromorphe non triviale de degr\'e $\geq 1$. Les sections
m\'eromorphes de $\F_{A_{u}}$ s'identifient aux solutions m\'eromorphes
$X = \begin{pmatrix} f \\ g \end{pmatrix}$ de l'\'equation $\sq X = A_{u} X$,
\ie\ du syst\`eme:
$$
\begin{cases}
\sq f = f + u g, \\
\sq g = z g.
\end{cases}
$$
Le degr\'e $\deg X \underset{def}{=} \deg \div_{\Eq} X$ de la section $X$ 
se calcule ainsi: si $h \in \Kw^{*}$ est telle que $h^{-1} X \in \Rw^{2}$ ne 
s'annule en aucun point, alors $\div_{\Eq} X \underset{def}{=} \div_{\Eq} h$
et $\deg \div_{\Eq} X = \deg \div_{\Eq} h$. Il est clair que, si $\phi$
est elliptique et non triviale, alors $\deg (\phi X) = \deg X$. \\
Si $g = 0$, $f$ est elliptique et, par hypoth\`ese, non triviale, 
donc $\deg X = \deg f = 0$ (ce cas correspond au sous-fibr\'e $\F_{(1)}$). \\
Si $g \neq 0$, alors $g = h \theta_{q}$, o\`u $h$ est elliptique et, quitte 
\`a diviser par $h$, on peut aussi bien supposer que $g = \theta_{q}$ et
$\sq f - f = u \theta_{q}$. Dire que $\deg X \geq 1$ c'est dire que $f$
et $g$ ont un z\'ero commun, donc que $f = \theta_{q} f_{1}$, o\`u
$f_{1} \in \Rw$. Mais alors $z \sq f_{1} - f_{1} = u$, qui n'a de 
solution dans $\Rw$ que si $u = 0$ (on le voit en examinant les s\'eries
de Laurent). \\
Nous avons donc une dichotomie: soit $u = 0$ et $A_{u}$ et $\F_{A_{u}}$
sont scind\'es; soit $u \neq 0$, et $A_{u}$ est ind\'ecomposable et
$\F_{A_{u}}$ est stable. Dans ce dernier cas, la filtration de Harder-
Narasimhan de $\F_{A_{u}}$ est triviale, alors que celle induite par
la filtration par les pentes de $A_{u}$ ne l'est pas.

\paragraph{Une famille semi-stable.}

Soit $A_{u} = \begin{pmatrix} 1 & u \\ 0 & z^{2} \end{pmatrix}$, 
avec $u \in \Ka$. D'apr\`es \cite{RSZ}, on peut ramener $A_{u}$ \`a 
la forme forme $A_{v}$ avec $v \in \C + \C z$ via un unique morphisme 
trivial sur les gradu\'es associ\'es. Nous supposerons donc d'embl\'ee
que $u = u_{0} + u_{1} z$. Par un raisonnement similaire \`a celui 
du paragraphe pr\'ec\'edent, on voit que la classe d'isomorphie de
$A_{u}$ d\'etermine $(u_{0},u_{1}) \in \C^{2}$ \`a un facteur pr\`es dans 
$\C^{*}$. De plus, si $u \neq 0$, l'objet $A_{u}$ est ind\'ecomposable. \\
Le degr\'e de $\F_{A_{u}}$ est $2$, sa pente est $1$. Si $u = 0$, ce
fibr\'e est scind\'e. On suppose d\'esormais que $u \neq 0$. On va voir
que, dans ce cas, $\F_{A_{u}}$ est semi-stable et non stable: autrement
dit, il admet des sections m\'eromorphes de degr\'e $1$, mais pas plus.
Une section m\'eromorphe de $\F_{A_{u}}$ est (par les identifications
habituelles) de la forme $X = \begin{pmatrix} f \\ g \end{pmatrix}$, o\`u:
$$
\begin{cases}
\sq f = f + u g, \\
\sq g = z^{2} g.
\end{cases}
$$
Si $g = 0$ et $f \neq 0$, $\deg X = 0$ (cas du sous-fibr\'e $\F_{(1)}$).
Si $g \neq 0$, on se ram\`ene au cas o\`u $g = \theta_{q}^{2}$, on
\'ecrit $f = f_{1} \theta_{q}^{2}$ et l'on doit avoir:
$z^{2} \sq f_{1} - f_{1} = u$. Pour que $\deg X \geq 1$, il faut
que $f_{1}$ ait au pire un p\^ole simple sur $\Eq$, soit $\overline{a}$,
donc que l'on ait $f_{1} = \dfrac{h}{\theta_{a}}$ avec $h \in \Rw$.
On est donc ramen\'e \`a chercher $a \in \C^{*}$ et $h \in \Rw$ tels
que $a z \sq h - h = u \theta_{a}$. Par d\'eveloppement en s\'erie
de Laurent $h = \sum h_{n} z^{n}$, cette \'equation \'equivaut \`a:
$$
\forall n \in \Z \;,\; a q^{n-1} h_{n-1} - h_{n} =
\left(u_{0} q^{-n(n+1)/2} + u_{1} a q^{-n(n-1)/2}\right) a^{-n},
$$
soit encore \`a:
$$
\forall n \in \Z \;,\; 
\dfrac{h_{n-1}}{q^{(n-1)(n-2)/2} a^{n-1}} - \dfrac{h_{n}}{q^{n(n-1)/2} a^{n}}
= \left(u_{0} q^{-n^{2}} + u_{1} a q^{-n(n-1)}\right) a^{-2n}.
$$
Il s'agit essentiellement d'une \emph{transformation de $q$-Borel}
(voir \cite{RSZ} et la section \ref{subsection:calculsdinvariants}).
Par sommation et annulation t\'elescopique, on
voit imm\'ediatement qu'une condition \emph{n\'ecessaire} est la
nullit\'e de:
$$
\phi_{u}(a) := 
\sum_{n \in \Z} \left(u_{0} q^{-n^{2}} + u_{1} a q^{-n(n-1)}\right) a^{-2n}.
$$
Les m\^emes majorations que dans \cite{JSSTO}, preuve du lemme 2.9,
montrent que c'est une condition suffisante. \\
On a sans peine:
$$
\phi(a) = 
u_{0} \theta_{q^{2}}(q a^{-2}) + u_{1} a^{-1} \theta_{q^{2}}(a^{-2}).
$$
Notons que la condition pos\'ee, $\phi(a) = 0$, est invariante par
$a \leftarrow q a$, comme il se doit (c'est une condition sur le p\^ole
$\overline{a} \in \Eq$). \\
Si $u_{0} = 0 \neq u_{1}$, on doit r\'esoudre $\theta_{q^{2}}(a^{-2}) = 0$, 
ce qui donne $a^{-2} \in [-1,q^{2}]$. Les deux solutions dans $\Eq$ sont 
$\overline{\ii}$ et $\overline{-\ii}$. De m\^eme, si $u_{0} \neq 0 = u_{1}$,
on r\'esoud $\theta_{q^{2}}(q a^{-2}) = 0$, donc $a^{-2} \in [-q,q^{2}]$,
donc les deux solutions dans $\Eq$ sont $\overline{\sqrt{-q}}$ et 
$\overline{-\sqrt{-q}}$ (avec un choix arbitraire de racine carr\'ee). \\
Supposons $u_{0} u_{1} \neq 0$. La fonction
$\psi(a) = \dfrac{\theta_{q^{2}}(q a^{-2})}{a^{-1} \theta_{q^{2}}(a^{-2})}$
est $q$-elliptique 
\footnote{En fait, on a une factorisation:
$\psi(a) = \dfrac{\theta_{q}(-\sqrt{-q} a) \, \theta_{q}(\sqrt{-q} a)}
{a \, \theta_{q}(-\ii a) \, \theta_{q}(\ii a)} \cdot$}
et admet, dans $\Eq$, deux p\^oles simples et deux 
z\'eros simples. Elle prend donc chaque valeur $\dfrac{u_{1}}{u_{0}}$
deux fois, et il y a encore deux p\^oles $\overline{a} \in \Eq$ possibles.\\
Si l'on a trouv\'e $\overline{a_{1}} \neq \overline{a_{2}}$, les sections
$X_{1}$ et $X_{2}$ correspondantes sont non proportionnelles. On peut en 
d\'eduire que $\F_{A_{u}}$ est scind\'e dans ce cas (qui est g\'en\'erique).



\section{Le ph\'enom\`ene de Stokes}
\label{section:phenomenedeStokes}


\subsection{Classification d'\'equations irr\'eguli\`eres}
\label{subsection:classification}

D'apr\`es la section \ref{subsubsection:equationspentesentieres},
on peut repr\'esenter tout objet $M$ de $\Ee$ par une matrice $A$
de la forme (\ref{equation:formestandard}) page 
\pageref{equation:formestandard}. 
Le gradu\'e $M_{0} = \gr M$
est alors d\'ecrit par la matrice:
\begin{equation} 
\label{equation:formestandardpure}
A_{0} = 
\begin{pmatrix}
z^{\mu_{1}} A_{1}  & \ldots & \ldots & \ldots & \ldots \\
\ldots & \ldots & \ldots  & 0 & \ldots \\
0      & \ldots & \ldots   & \ldots & \ldots \\
\ldots & 0 & \ldots  & \ldots & \ldots \\
0      & \ldots & 0       & \ldots & z^{\mu_{k}} A_{k}    
\end{pmatrix},
\end{equation}
De la derni\`ere assertion du th\'eor\`eme \ref{thm:filtration},
il d\'ecoule alors qu'il existe une unique matrice \emph{formelle},
\cad\ \`a coefficients dans $\Kf$:
\begin{equation} 
\label{equation:automorphisme}
F = 
\begin{pmatrix}
I_{r_{1}} & \ldots & \ldots & \ldots & \ldots \\
\ldots & \ldots & \ldots  & F_{i,j} & \ldots \\
0      & \ldots & \ldots   & \ldots & \ldots \\
\ldots & 0 & \ldots  & \ldots & \ldots \\
0      & \ldots & 0       & \ldots & I_{r_{k}}     
\end{pmatrix},
\end{equation}
telle que $F[A_{0}] = A$. Nous la noterons $\hat{F}_{A} \in \G(\Kf)$,
la notation $\G$ d\'esignant le sous-groupe alg\'ebrique unipotent
de $\GL_{n}$ correspondant au format ci-dessus. Plus g\'en\'eralement,
si $A'$ est aussi de la forme (\ref{equation:formestandard}), alors
l'unique $F \in \G(\Kf)$ tel que $F[A] = A'$ est
$\hat{F}_{A,A'} = \hat{F}_{A'} \left(\hat{F}_{A}\right)^{-1}$.
Les isomorphismes formels $\hat{F}_{A}$, $\hat{F}_{A,A'}$ sont
en g\'en\'eral tr\`es divergents: les blocs $F_{i,j}$ tels que
$\mu_{j} - \mu_{i} = \delta \in \N^{*}$ sont en g\'en\'eral de
\emph{niveau $q$-Gevrey $\delta$}, autrement dit, leurs coefficients
sont de la forme $\sum f_{n} z^{n}$ avec 
$f_{n} = O\left(R^{n} q^{n^{2}/2 \delta}\right)$ pour un certain $R > 0$.

\paragraph{Classification analytique isoformelle.}

Il s'agit de la classification analytique \`a classe formelle fix\'ee
(ou \`a gradu\'e fix\'e, ce qui revient au m\^eme). On fixe un
module pur $M_{0} = P_{1} \oplus \cdots \oplus P_{k}$ dans $\Epe$
et l'on consid\`ere tous les couples $(M,u)$, o\`u $M$ est un objet 
de $\Ee$ et $u: \gr M \rightarrow M_{0}$ un isomorphisme. On dit que 
$(M,u) \sim (N,v)$ s'il existe un isomorphisme $\phi: M \rightarrow N$
tel que $v \circ (\gr \phi) = u$ (isomorphismes triviaux sur le gradu\'e).
Avec les notations vues plus haut, il n'est pas tr\`es difficile de
voir que, si $M_{0}$ est d\'ecrit par la matrice $A_{0}$ de la forme
(\ref{equation:formestandardpure}), alors un couple $(M,u)$ est la
m\^eme chose qu'une matrice $A$ sous la forme (\ref{equation:formestandard})
et que, si $(M,u)$ et $(M',u')$ correspondent \`a $A$ et $A'$, alors
ils sont \'equivalents si, et seulement si $\hat{F}_{A,A'} \in \G(\Ka)$. \\

Il a \'et\'e d\'emontr\'e dans \cite{RSZ} que l'ensemble $\F(M_{0})$
des classes pour cette relation est un espace affine de dimension
$\irr(M_{0}) \underset{def}{=} 
\sum\limits_{1 \leq i < j \leq k} r_{i} r_{j} (\mu_{j} - \mu_{i})$,
o\`u les $\mu_{i}$, $r_{i}$ proviennent du polygone de Newton de $M_{0}$
(et de tous les $M$ de sa classe formelle)
\footnote{Ce r\'esultat est utilis\'e dans \cite{vdP}, mais attribu\'e
de mani\`ere erronn\'ee \`a Ramis et Sauloy, et donn\'e sans r\'ef\'erence.}. 
Il s'agit en fait d'un vrai sch\'ema de modules pour ce probl\`eme.
Une preuve consiste \`a v\'erifier que, dans chaque classe, il existe
une unique matrice $A$ en \emph{forme normale}, \cad\ en forme standard
(\ref{equation:formestandard}) et telle que chaque $U_{i,j}$ est \`a
coefficients dans $\sum\limits_{\mu_{i} \leq d < \mu_{j}} \C z^{d}$.
Cette forme normale est inspir\'ee de Birkhoff et Guenther.
Cependant, l'objectif de \cite{RSZ} est d'obtenir des \emph{invariants
transcendants} sous une forme $q$-analogue aux th\'eor\`emes de
Malgrange-Sibuya. Nous allons d\'ecrire l'un de ces r\'esultats
sous la forme moins puissante, mais plus simple, de \cite{JSSTO}.

\paragraph{Sommation de $\hat{F}_{A}$.}

Dans \cite{JSSTO}, on d\'efinit un sous-ensemble fini explicite
$\Sigma_{A_{0}}$ de $\Eq$, qui est, en g\'en\'eral, de cardinal
$\irr(A_{0})$ (et moins dans certains cas ``r\'esonnants'').
On prouve alors:

\begin{thm}
Pour tout $\overline{c} \in \Eq \setminus \Sigma_{A_{0}}$, il existe
un unique isomorphisme m\'eromorphe $F: A_{0} \rightarrow A$ tel que 
$F \in \G(\Kw)$, dont les p\^oles sont situ\'es sur la $q$-spirale
discr\`ete $[-c;q] = -c q^{\Z}$ et tel que, pour $1 \leq i < j \leq k$, 
les p\^oles de $F_{i,j}$ ont une multiplicit\'e $\leq \mu_{j} - \mu_{i}$.
(Ce que l'on peut \'ecrire: 
$\div_{\Eq} F_{i,j} \geq - (\mu_{j} - \mu_{i}) [\overline{c}]$.)
Nous noterons $S_{\overline{c}} \hat{F}_{A}$ cette matrice $F$.
\end{thm}

La description de $\Sigma_{A_{0}}$ et le calcul de 
$S_{\overline{c}} \hat{F}_{A}$ seront explicit\'es dans un cas
crucial \`a la section \ref{subsection:calculsdinvariants}, en
vue de comparaison avec d'autres constructions int\'eressantes. \\

\Rem
Le morphisme m\'eromorphe $S_{\overline{c}} \hat{F}_{A}$ induit
un isomorphisme holomorphe du faisceau $\F_{A_{0}}$ sur 
le faisceau $\F_{A}$ en restriction \`a l'ouvert 
$\Eq \setminus \{\overline{c}\}$. ces deux faisceaux sont 
donc localement isomorphes. Comme il est facile de prouver
que $\F_{A_{0}}$ est un fibr\'e, cela fournit une nouvelle
preuve du fait que $\F_{A}$ l'est \'egalement. En fait, cela
montre que le fibr\'e $\F_{A}$ s'obtient \`a partir du fibr\'e
$\F_{A_{0}}$ par l'op\'eration de cohomologie non-ab\'elienne
``torsion par un cocycle'' (\cite{Frenkel}). \\

Nous consid\'erons $S_{\overline{c}} \hat{F}_{A}$ comme 
une \emph{sommation de la s\'erie divergente $\hat{F}_{A}$ 
dans la direction $\overline{c} \in \Eq$}. Ce point de vue admet
diverses justifications, dont celle-ci: selon la th\'eorie asymptotique
$q$-Gevrey d\'evelopp\'ee par Ramis et Zhang dans \cite{RZ}, 
$S_{\overline{c}} \hat{F}_{A}$ est asymptote \`a $\hat{F}_{A}$. \\

Il est facile de d\'eduire du th\'eor\`eme que, si $A$ et $A'$
sont sous la forme (\ref{equation:formestandard}) (avec m\^eme
gradu\'e $A_{0}$), alors 
$S_{\overline{c}} \hat{F}_{A'} 
\left(S_{\overline{c}} \hat{F}_{A}\right)^{-1}$
est l'unique isomorphisme m\'eromorphe de $A$ dans $A'$
satisfaisant aux m\^emes conditions de polarit\'e: il est donc
l\'egitime de le noter $S_{\overline{c}} \hat{F}_{A,A'}$. On peut
alors d\'emontrer:

\begin{prop}
Soit $\overline{c} \in \Eq \setminus \Sigma_{A_{0}}$. Pour que $A$ 
et $A'$ soient dans la m\^eme classe analytique, il faut, et il suffit, 
que $S_{\overline{c}} \hat{F}_{A,A'}$ n'ait pas de p\^ole sur $[-c;q]$.
\end{prop}

Naturellement, dans ce cas, $\hat{F}_{A,A'}$ est analytique et tous 
les $S_{\overline{c}} \hat{F}_{A,A'}$ lui sont \'egaux. Ce qui est
remarquable dans ce r\'esultat, c'est que l'absence des p\^oles sur
$[-c;q]$, qui \emph{a priori} ne devrait entrainer que l'holomorphie
sur $\C^{*}$, suffit en fait \`a garantir la m\'eromorphie en $0$.
L'argument est directement li\'e \`a celui qui a permis de prouver
le th\'eor\`eme \ref{thm:pleinefidelite}.

\paragraph{Op\'erateurs de Stokes.}

Classiquement, l'ambigu\"it\'e dans la sommation d'une solution
divergente d'\'equation fonctionnelle est le \emph{ph\'enom\`ene
de Stokes}. Il se r\'ealise ici sous la forme suivante. On prend
deux directions de sommation autoris\'ees 
$\overline{c}, \overline{d} \in \Eq \setminus \Sigma_{A_{0}}$
et l'on pose:
$$
S_{\overline{c},\overline{d}} \hat{F}_{A} =
\left(S_{\overline{c}} \hat{F}_{A}\right)^{-1} S_{\overline{d}} \hat{F}_{A}.
$$
C'est donc un automorphisme m\'eromorphe de $A_{0}$. Il est holomorphe
sur l'ouvert $U_{c,d} = \C^{*} \setminus [-c,-d;q]$
de $\C^{*}$. On peut aussi bien le consid\'erer comme un automorphisme
m\'eromorphe du fibr\'e $\F_{A_{0}}$, holomorphe sur l'ouvert 
$V_{\overline{c}, \overline{d}} = 
\Eq \setminus \{\overline{c}, \overline{d}\}$ de $\Eq$.
On a donc $U_{c,d} = \pi^{-1}\left(V_{\overline{c}, \overline{d}}\right)$. \\

On d\'efinit le faisceau en groupes (non commutatifs) $\Lambda_{I}(M_{0})$
comme suit; pour tout ouvert $V$ de $\Eq$:
$$
\Gamma\left(V,\Lambda_{I}(M_{0})\right) =
\{F \in \G\left(\O(\pi^{-1}(V))\right) \mid F[A_{0}] = A_{0}\}.
$$
C'est le faisceau des \emph{automorphismes de $M_{0}$ tangents 
\`a l'identit\'e}. La terminologie (emprunt\'ee au cas analogue
des \'equations diff\'erentielles) est justifi\'ee par le fait
que, si $F \in \Gamma\left(V,\Lambda_{I}(M_{0})\right)$, alors
$F - I_{n}$ est plat pr\`es de $0$ sur $V$. En effet, tout bloc
hors diagonal $F_{i,j}$ v\'erifie 
$\sq F_{i,j} = z^{\mu_{i} - \mu_{j}} A_{i} F_{i,j} A_{j}^{-1}$,
d'o\`u l'on tire que $\theta_{q}^{\mu_{j} - \mu_{i}} F_{i,j}$
est \`a croissance mod\'er\'ee pr\`es de $0$ sur son ouvert de
d\'efinition. On dit que $F_{i,j}$ est \emph{$t$-plat}, o\`u
$t = \mu_{j} - \mu_{i}$. \\

On a a donc:
$$
S_{\overline{c},\overline{d}} \hat{F}_{A} \in
\Gamma\left(V_{\overline{c}, \overline{d}},\Lambda_{I}(M_{0})\right).
$$
Notant $V_{\overline{c}} = \Eq \setminus \{\overline{c}\}$, 
il est donc clair que l'on a d\'efini un \'el\'ement de l'ensemble
$Z^{1}\left(\V,\Lambda_{I}(M_{0})\right)$ des cocycles du faisceau
en groupes $\Lambda_{I}(M_{0})$ associ\'es au recouvrement 
$\V = (V_{\overline{c}})$ de $\Eq$ 
($Z^{1}\left(\V,\Lambda_{I}(M_{0})\right)$
est bien un ensemble point\'e, pas un groupe). 
Voici le th\'eor\`eme de type Malgrange-Sibuya annonc\'e:

\begin{thm} 
\label{thm:q-Malgrange-Sibuya}
On obtient ainsi une bijection:
$$
\F(M_{0}) \simeq H^{1}\left(\Eq,\Lambda_{I}(M_{0})\right).
$$
\end{thm}

\paragraph{D\'evissage $q$-Gevrey.}

On peut plagier le \emph{d\'evissage Gevrey} de Ramis en un
d\'evissage $q$-Gevrey du faisceau $\Lambda_{I}(M_{0})$ et
de l'ensemble de cohomologie $H^{1}\left(\Eq,\Lambda_{I}(M_{0})\right)$.
Tout d'abord, il s'agit d'un faisceau de groupes unipotents
dont on peut facilement d\'ecrire le faisceau $\lambda_{I}$
des alg\`ebres de Lie:
$$
\Gamma\left(V,\lambda_{I}(M_{0})\right) =
\{F \in \g\left(\O(\pi^{-1}(V))\right) \mid (\sq F) A_{0} = A_{0} F\}.
$$
Naturellement, $\g$ d\'esigne l'alg\`ebre de Lie de $\G$, d\'efinie par:
$$
F \in \g(K) \Leftrightarrow I_{n} + F \in \G(K)
\Leftrightarrow \exp(F) \in \G(K).
$$
On note $\lambda_{I}^{t}(M_{0})$ le sous-faisceau en alg\`ebres
de Lie nilpotentes de $\lambda_{I}(M_{0})$ form\'e des \'el\'ements
$t$-plats. On peut montrer qu'il s'agit des \'el\'ements dont
les seuls blocs $F_{i,j}$ non nuls sont ceux tels que 
$ \mu_{j} - \mu_{i} \geq t$ (donc triangulaires sup\'erieurs 
par blocs, ``assez loin'' de la diagonale). \\
Si l'on note de m\^eme $\lambda_{I}^{(t)}(M_{0})$ le sous-faisceau
form\'e des \'el\'ements dont les seuls blocs $F_{i,j}$ non nuls 
sont ceux tels que $ \mu_{j} - \mu_{i} = t$ (donc ne comportant
qu'une ``surdiagonale par blocs''), on constate que $\lambda_{I}(M_{0})$
est gradu\'ee:
$$
\lambda_{I}(M_{0}) = \bigoplus_{t \in \N^{*}} \lambda_{I}^{(t)}(M_{0}).
$$
Les $\lambda_{I}^{t}(M_{0})$ forment la filtration correspondante
par des id\'eaux:
$$
\lambda_{I}^{t}(M_{0}) = \bigoplus_{t' \geq t} \lambda_{I}^{(t')}(M_{0}).
$$

\begin{prop}
(i) Le faisceau $\lambda_{I}^{(t)}(M_{0})$ est localement libre.
C'est le fibr\'e associ\'e au module pur isocline:
$\underline{End}(M_{0})_{(-t)}$. \\
(ii) Le faisceau $\lambda_{I}(M_{0})$ est localement libre.
C'est le fibr\'e associ\'e au module pur: $\underline{End}(M_{0})_{<0}$.
\end{prop}

Rappelons que $\underline{Hom}$ d\'esigne le \emph{Hom interne},
$\underline{End}(M_{0})$ n'est autre que le module 
$\underline{Hom}(M_{0},M_{0}) = \bigoplus \underline{Hom}(P_{i},P_{j})$
(o\`u $M_{0} = \bigoplus P_{i}$, chaque $P_{i}$ \'etant pur isocline
de pente $\mu_{i}$). On a alors:
$\underline{End}(M_{0})_{(t)} = 
\bigoplus\limits_{\mu_{j} - \mu_{i} = t} \underline{Hom}(P_{i},P_{j})$. \\

\`A la \emph{graduation} de $\lambda_{I}(M_{0})$ correspond 
une \emph{filtration} de $\Lambda_{I}(M_{0})$ par les (faisceaux en)
sous-groupes distingu\'es:
$$
\Lambda_{I}^{t}(M_{0}) = I_{n} + \lambda_{I}^{t}(M_{0}) =
\exp \lambda_{I}^{t}(M_{0}).
$$
Le \emph{d\'evissage $q$-Gevrey} du faisceau en groupes unipotents
non commutatifs par des fibr\'es est donn\'e par les suites exactes:
\begin{equation} 
\label{equation:suite-exacte}
1 \rightarrow \Lambda_{I}^{t+1}(M_{0}) \rightarrow
\Lambda_{I}^{t}(M_{0}) \rightarrow \lambda_{I}^{(t)}(M_{0})
\rightarrow 0.
\end{equation}
Il y a aussi une suite d'extensions centrales:
\begin{equation} \label{equation:extension-centrale}
0 \rightarrow \lambda_{I}^{(t)}(M_{0}) \rightarrow
\dfrac{\Lambda_{I}(M_{0})}{\Lambda_{I}^{t+1}(M_{0})} \rightarrow
\dfrac{\Lambda_{I}(M_{0})}{\Lambda_{I}^{t}(M_{0})} \rightarrow 1.
\end{equation}
Ce sont ces suites qui permettent une preuve ``\'el\'ementaire''
du th\'eor\`eme \ref{thm:q-Malgrange-Sibuya}.


\subsection{Calculs explicites d'invariants}
\label{subsection:calculsdinvariants}

Puisque l'espace de modules $\F(M_{0})$ est de dimension finie, on doit
pouvoir d\'eterminer des coordonn\'ees, \ie\ un jeu complet d'invariants
finis. Les $\irr(M_{0})$ coefficients des polyn\^omes $U_{i,j}$ dans
l'\'ecriture en forme normale conviennent, mais ne sont pas tr\`es
int\'eressants. \\

On va d\'etailler ici le cas d'un polygone de Newton \`a deux pentes.
Trois types d'invariants se pr\'esentent. Le premier est fourni par
la transformation de $q$-Borel (qui, en analyse, devrait \^etre
accompagn\'ee de la ransformation de $q$-Laplace, \cite{Zha02}).
En th\'eorie de Galois apparaissent les $q$-d\'eriv\'ees \'etrang\`eres
(\cite{RS1} et \cite{RS2}). Enfin, nous allons rencontrer une application
inattendue de la dualit\'e de Serre. Actuellement, nous ne savons
d\'efinir pour un nombre arbitraire de pentes que les $q$-d\'eriv\'ees 
\'etrang\`eres (et les pentes doivent \^etre enti\`eres).

\begin{lemma}
Pour tout module $M$, on a une identification naturelle:
$$
\Gamma^{1}(M) \simeq H^{1}(\Eq,\F_{M}).
$$
\end{lemma}
\Pr
Le foncteur exact \`a gauche $M \leadsto \Gamma(M)$ s'identifie 
naturellement au foncteur compos\'e du foncteur exact $M \leadsto \F_{M}$
et du foncteur exact \`a gauche $\F \leadsto \Gamma(\Eq,\F)$. Leurs
foncteurs d\'eriv\'es \`a droite sont donc \'egaux. 
\hfill $\Box$ \\

Soit $M_{0} = P_{1} \oplus P_{2}$, o\`u $P_{i} = M_{z^{\mu_{i}} A_{i}}$,
$A_{i} \in \GL_{r_{i}}(\C)$ ($i = 1,2$), $\mu_{1} < \mu_{2} \in \Z$.
L'espace $\F(M_{0})$ s'identifie naturellement \`a $Ext(P_{2},P_{1})$
et ce dernier, d'apr\`es la section \ref{subsubsection:modules}, \`a
$\Gamma^{1}(P)$, o\`u $P = P_{2}^{\vee} \otimes P_{1}$ est un module
pur de pente $\mu = \mu_{1} - \mu_{2} < 0$ et de rang $r = r_{1} r_{2}$.
Pour l'\'etude du cas de deux pentes, on peut donc se restreindre aux
calculs sur $M_{0} = P \oplus \underline{1} = M_{A_{0}}$, o\`u
$P = M_{z^{\mu}} A$:
$$
A_{0} = \begin{pmatrix} z^{\mu} A & 0 \\ 0 & 1 \end{pmatrix}, \quad
\mu \in \Z, \mu <0 \text{~et~} A \in \GL_{r}(\C), r \in \N^{*}.
$$
Nous poserons $d = - \mu$ (``niveau'' $q$-Gevrey). Nous allons 
expliciter dans ce cas les isomorphismes 
$\F(M_{0}) = \Gamma^{1}(P) \simeq H^{1}(\Eq,\F_{M_{0}})$. \\

Les sections de $\Lambda_{I}(M_{0})$ sont ici les matrices 
$\begin{pmatrix} I_{r} & F \\ 0 & 1 \end{pmatrix}$ telles que
$\sq F = z^{\mu} A F$, la multiplication dans le groupe 
$\Gamma\left(V,\Lambda_{I}(M_{0})\right)$ correspondant \`a l'addition 
de leurs composantes $F$. Le faisceau $\Lambda_{I}(M_{0})$ est donc
isomorphe au fibr\'e $\F_{P} = \F_{z^{\mu} A}$. Puisque l'on a
un faisceau en groupes commutatifs, le $H^{1}$ est ici un groupe 
(cohomologie des faisceaux ab\'eliens).

\paragraph{Calcul de cocycles.}

On va calculer des cocycles par ``sommation''. Pour cela, on choisit 
un module $M_{U}$, de  matrice 
$A_{U} = \begin{pmatrix} z^{\mu} A & U \\ 0 & 1 \end{pmatrix}$
de gradu\'e $A_{0}$, avec $U \in \Ma_{r,1}(\Ka)$; on impose
que $U$ soit polynomial (en fait, holomorphe sur $\C^{*}$ et
m\'eromorphe en $0$ suffit pour les calculs qui suivent).
On fixe une direction de sommation $\overline{c} \in \Eq$,
et l'on cherche $\begin{pmatrix} I_{r} & F \\ 0 & 1 \end{pmatrix}$ 
qui soit un isomorphisme m\'eromorphe de $A_{0}$ dans $A_{U}$, avec 
pour seuls p\^oles $[-c;q]$, de multiplicit\'e $\geq d = - \mu$.
L'\'equation satisfaite par $F$ est:
$$
\sq F - z^{\mu} A F = U.
$$
On pose $F = \dfrac{G_{c}}{\theta_{q,c}^{d}}$. On est ramen\'e \`a
r\'esoudre l'\'equation:
$$
c^{d} \sq G_{c} - A G_{c} = U \theta_{q,c}^{d}.
$$
On d\'eveloppe en s\'eries de Laurent $G_{c} = \sum G_{n} z^{n}$,
$V = U \theta_{q,c}^{d} = \sum V_{n} z^{n}$, d'o\`u:
$$
\forall n \in \Z \;,\; (c^{d} q^{n} I_{r} - A) G_{n} = V_{n}.
$$
Supposons maintenant que toutes les matrices $c^{d} q^{n} I_{r} - A$
soient inversibles, \ie\ \'el\'ements de $\GL_{r}(\C)$. Cela \'equivaut \`a:
$$
\overline{c} \in \Eq \setminus \Sigma_{A_{0}}, \quad \text{~o\`u~}
\Sigma_{A_{0}} = \overline{\sqrt[d]{\Sp(A)}} = 
\bigl(\sqrt[d]{\Sp(A)} \pmod{q^{\Z}}\bigr).
$$
Notons que, g\'en\'eriquement (\cad\ si $\Sp(A)$ admet $r$ valeurs
distinctes modulo $q^{\Z}$), l'ensemble interdit $\Sigma_{A_{0}}$
a $rd = \irr(M_{0})$ \'el\'ements. Pour une direction autoris\'ee,
les \'equations ci-dessus admettent une unique solution:
\begin{eqnarray*}
G_{c} & = & \sum_{n \in \Z} (c^{d} q^{n} I_{r} - A)^{-1} V_{n} z^{n}, \\
F = F_{\overline{c}} & = & \dfrac{1}{\theta_{q,c}^{d}} 
\sum_{n \in \Z} (c^{d} q^{n} I_{r} - A)^{-1} V_{n} z^{n}.
\end{eqnarray*}
La notation $F_{\overline{c}}$ est l\'egitime, car cette fonction
ne d\'epend que de la classe de $c$ dans $\Eq$. Si l'on pose:
$$
F_{\overline{c},\overline{d}} = F_{\overline{d}} - F_{\overline{c}},
$$
on voit que $F_{\overline{c},\overline{d}}$ est holomorphe sur 
$\C^{*} \setminus [-c,-d;q]$ et v\'erifie
$\sq F_{\overline{c},\overline{d}} = z^{\mu} A F_{\overline{c},\overline{d}}$.
Les fonctions $F_{\overline{c},\overline{d}}$ sont des sections du fibr\'e 
$\F_{P} = \F_{z^{\mu} A}$ sur les ouverts $V_{\overline{c},\overline{d}}$,
et constituent un cocycle. La classe de ce cocycle dans le groupe de 
cohomologie $H^{1}(\Eq,\F_{M_{0}})$ s'identifie \`a la fois \`a la classe 
de $A_{U}$ dans $\F(M_{0})$ et \`a la classe de l'extension $M_{U}$ de 
$\underline{1}$ par $P$ dans $\Gamma^{1}(P) = Ext(\underline{1},P)$.
Nous noterons $\cl(M_{U})$ cette classe (dans l'un quelconque des
trois ensembles ainsi identifi\'es). 

\paragraph{Invariants \`a la $q$-Borel.}

Pour r\'esoudre (ou pour \'etudier l'obtruction \`a r\'esoudre)
l'\'equation $\sq F - z^{\mu} A F = U$, on d\'eveloppe en s\'eries
de Laurent: $F = \sum F_{n} z^{n}$, $U = \sum U_{n} z^{n}$, d'o\`u:
$$
\forall n \in \Z \;,\; q^{n} F_{n} - A F_{n+d} = U_{n}.
$$
On introduit des coefficients $t_{n}$ tels que:
$$
\forall n \in \Z \;,\; q^{n} t_{n} = t_{n-d}.
$$
Par exemple, le d\'eveloppement en s\'erie de Laurent 
$\theta_{q}^{d} = \sum t_{n} z^{n}$ fournit de tels coefficients,
comme il d\'ecoule de l'\'equation fonctionnelle 
$\sq \theta_{q}^{d} = z^{d} \theta_{q}^{d}$.
On appelle \emph{transform\'ees de $q$-Borel au niveau $d$}
d'une s\'erie $f(z) = \sum\limits_{n \in \Z} f_{n} z^{n}$ la s\'erie:
$$
\Bd f(\xi) = \sum_{n \in \Z} q^{-n} t_{-n} f_{n} \xi^{n}.
$$
Un petit calcul montre que notre relation ci-dessus \'equivaut \`a:
$$
(I_{r} - \xi^{-d} A) \Bd F(q \xi) = \Bd U(\xi).
$$
Si $U$ (resp. $F$) est analytique, $\Bd U$ (resp. $\Bd F$)
est enti\`ere, et l'on peut \'evaluer cette \'egalit\'e
partout. On choisit une matrice $B$ racine $d$\ieme\ de $A$,
et l'on voit qu'une condition n\'ecessaire est l'annulation
des $\Bd(j \B)$ pour $j^{d} = 1$. C'est en fait une condition
suffisante, et l'on peut d\'emontrer que les $d$ vecteurs
$\Bd(j \B) \in \C^{r}$ constituent un jeu complet d'invariants
analytiques. Plus pr\'ecis\'ement, l'application qui, \`a $U$, 
associe le $d$-uplet des $\Bd(j \B)$ induit un isomorphisme de 
l'espace vectoriel $\F(M_{0})$ sur $\C^{rd}$.

\paragraph{Invariants de \cite{RS1}, \cite{RS2}.}

On fixe un $a \in \C^{*}$ ``g\'en\'erique'' (pratiquement, loin
de tous les points qui vont intervenir). Il jouera le r\^ole
d'un point-base. On fixe une direction de sommation autoris\'ee
arbitraire $\overline{c_{0}} \in \Eq$. L'application
$\overline{d} \mapsto F_{\overline{c_{0}},\overline{d}}(a)$
est m\'eromorphe sur $\Eq$, avec des p\^oles sur $\Sigma_{A_{0}}$. 
Elle est \`a valeurs dans \emph{l'alg\`ebre de Lie $\st(M_{U})$
du groupe de Stokes $\St(M_{U})$}. (C'est pour cela que l'on a
d\^u introduire $\overline{c_{0}}$, qui n'intervient pas dans le
r\'esultat du calcul). La prise de r\'esidu est une int\'egration, 
donc donne un r\'esultat dans l'espace vectoriel $\st(M_{U})$.
On pose, pour tout $\overline{c} \in \Sigma_{A_{0}}$:
$$
\Derc(A_{U}) = 
\Res_{\overline{d} = \overline{c}} F_{\overline{c_{0}},\overline{d}}(a).
$$
Les $\Derc(A_{U})$ prennent leurs valeurs dans des espaces vectoriels
de dimension totale $rd$ et constituent un jeu complet d'invariants
analytiques. Dans \cite{RS1}, on donne des formules de transformations
entre ces ``$q$-d\'eriv\'ees \'etrang\`eres'' et les invariants de
$q$-Borel.

\paragraph{Invariants provenant de la dualit\'e de Serre.}

Il s'agit de la dualit\'e de Serre pour les fibr\'es vectoriels
holomorphes sur une surface de Riemann compacte (\cite{Gunning}).
Le diviseur canonique de $\Eq$ est trivial. Il y a en effet une
diff\'erentielle globale de diviseur nul: c'est, par exemple
$\dfrac{dz}{z} = 2 \ii \pi dx$, o\`u $x$ et $z$ sont les uniformisantes
globales provenant des rev\^etements $\C \rightarrow \Eq$ et 
$\C^{*} \rightarrow \Eq$. De la dualit\'e de Serre, on d\'eduit
alors que, pour tout fibr\'e sur $\Eq$, les espaces vectoriels 
$H^{0}(\Eq,\F^{\vee})$ et $H^{1}(\Eq,\F)$ sont duaux l'un de l'autre.
Une base de $H^{0}(\Eq,\F_{M_{0}}^{\vee})$ fournira donc un syst\`eme
de coordonn\'ees sur $H^{1}(\Eq,\F_{M_{0}})$, \cad\ un jeu complet
d'invariants analytiques pour $\F(M_{0})$. \\

Le m\'ecanisme de dualit\'e est le suivant. Soit $X$ une section
globale de $\F_{M_{0}}^{\vee}$: c'est donc une solution holomorphe
sur $\C^{*}$ de l'\'equation:
$$
\sq X = z^{-\mu} \; {}^{t} A^{-1} X.
$$
Alors $Y = {}^{t} X$ est holomorphe sur $\C^{*}$ et v\'erifie:
$\sq Y = z^{-\mu} \, Y A^{-1}$. On forme le produit scalaire
$$
\phi_{\overline{c},\overline{d}} =
{<} X, F_{\overline{c},\overline{d}} {>} = 
{}^{t} X F_{\overline{c},\overline{d}} =
Y F_{\overline{c},\overline{d}}.
$$
C'est une fonction m\'eromorphe sur $\C^{*}$ avec des p\^oles connus, et, 
des \'equations satisfaites par $Y$ et $F_{\overline{c},\overline{d}}$,
on d\'eduit:
$$
\sq\phi_{\overline{c},\overline{d}}  = 
\sq Y \sq F_{\overline{c},\overline{d}} =
z^{-\mu} \, Y A^{-1} \, z^{\mu} A F_{\overline{c},\overline{d}} = 
\phi_{\overline{c},\overline{d}}.
$$
On a donc une fonction elliptique $\phi_{\overline{c},\overline{d}}$ 
de p\^oles $\overline{c},\overline{d} \in \Eq$. La somme de ses 
r\'esidus sur $\Eq$ est donc nulle. Par ailleurs, cette fonction 
est, par d\'efinition, une diff\'erence:
$$
\phi_{\overline{c},\overline{d}} = 
{<} X, F_{\overline{d}} - F_{\overline{c}} {>} =
{<} X, F_{\overline{d}} {>} - {<} X, F_{\overline{c}} {>} =
\phi_{\overline{d}} - \phi_{\overline{c}},
$$
o\`u chaque section $\phi_{\overline{c}} = {<} X, F_{\overline{c}} {>}$
a un seul p\^ole sur $\Eq$ (une $q$-spirale sur $\C^{*}$), \`a savoir
$\overline{-c}$. On a donc:
$$
\Res_{\overline{-c}} \, \phi_{\overline{c}} = 
\Res_{\overline{-d}} \, \phi_{\overline{d}}.
$$
(On prendra garde qu'il s'agit ici de r\'esidus par rapport
\`a la variable $z$ et non par rapport \`a la direction de
sommation, comme dans le cas des $\Derc$ !)
Par d\'efinition, le nombre calcul\'e ci-dessus, qui ne d\'epend pas
de $\overline{c} \in \Eq$, est celui associ\'e par la dualit\'e de Serre
\`a $X \in H^{0}(\Eq,\F_{M_{0}}^{\vee})$ et \`a 
$\cl(M_{U}) \in H^{1}(\Eq,\F_{M_{0}})$.
Nous le noterons ${<} X , \cl(M_{U}) {>}$.

\begin{lemma}
On a les \'egalit\'es:
$$
{<} X , \cl(M_{U}) {>} = [(\sq X) U]_{0} = 
[z^{d} {}^{t} X A^{-1} U]_{0} = [{}^{t} X A^{-1} U]_{\mu}.
$$
\end{lemma}
\Pr
Nous employons la notation commode suivante: si $f = \sum f_{n} z^{n}$,
alors $[f]_{n} = f_{n}$. Ici, $X = \sum X_{n} z^{n}$, 
$U = \sum U_{n} z^{n}$ (s\'eries de Laurent convergentes sur $\C^{*}$)
et l'on calcule le coefficient de degr\'e $0$ (resp. de degr\'e $\mu$)
de $z^{d} {}^{t} X A^{-1} U$ (resp. de ${}^{t} X A^{-1} U$). La derni\`ere
\'egalit\'e est donc triviale. La seconde l'est \'egalement, puisque
$\sq X = z^{d} \; {}^{t} X A^{-1}$. \\
Pour calculer le r\'esidu en $\overline{-c} \in \Eq$,
on choisit un repr\'esentant $c \in \C^{*}$ et une couronne 
fondamentale contenant le p\^ole $-c$. On \'ecrit la fronti\`ere
orient\'ee de cette couronne sous la forme $q \gamma - \gamma$, o\`u 
$\gamma$ est un cercle de centre $0$ parcouru positivement. Alors,
par le th\'eor\`eme des r\'esidus:
$$
{<} X , \cl(M_{U}) {>} = \Res_{\overline{-c}} \, \phi_{\overline{c}} =
\int_{q \gamma} {}^{t} X F_{\overline{c}} \dfrac{dz}{2 \ii \pi z}
- \int_{\gamma} {}^{t} X F_{\overline{c}} \dfrac{dz}{2 \ii \pi z} =
\int_{\gamma} 
\left(\sq({}^{t} X F_{\overline{c}}) - {}^{t} X F_{\overline{c}}\right)
\dfrac{dz}{2 \ii \pi z} = 
\int_{\gamma} (\sq X) U  \dfrac{dz}{2 \ii \pi z} = [(\sq X) U]_{0}.
$$
\hfill $\Box$ \\

Donc ${<} X , \cl(M_{U}) {>} = \sum q^{-n} X_{n} U_{n}$.
Pour produire concr\`etement de tels nombres, il est naturel
de construire la section $X$ \`a l'aide de fonctions Theta.
On choisit une racine $d$\ieme\ $B$ de la matrice $A$, et l'on 
voit que: 
$$
T_{B} = \Theta^{d}_{B}(z) \underset{def}{=} 
\theta^{d}(B^{-1} z) = \sum_{n \in \Z} t_{n} B^{-n} z^{n}
$$
v\'erifie $\sq T_{B} = T_{B} z^{d} A^{-1}$. On peut donc prendre
pour ${}^{t} X$ l'une quelconque des $r$ lignes de $T_{B}$.
En faisant le m\^eme calcul pour chacune des $d$ matrices $j B$
($j^{d} = 1$), on trouve une base de $H^{0}(\Eq,\F_{M_{0}}^{\vee})$
et les invariants correspondants sont les invariants de $q$-Borel
calcul\'es pr\'ec\'edemment.


\section{Constructions globales}
\label{section:constructionsglobales}

Historiquement, les \'equations aux $q$-diff\'erences int\'eressantes
sont d\'efinies sur la sph\`ere de Riemann $\Sr$, \cad\ \`a coefficients
rationnels: matrice $A \in \GL_{n}(\Kr)$ codant un objet $M_{A}$ de
$\Er = \DMr$ (voir \cite{Birkhoff1}). Par exemple, l'\'equation 
$q$-hyperg\'eom\'etrique (\cite{GR}, \cite{JSAIF}, \cite{DRSZ}). \\

Pour \'etudier une \'equation globale, on la localise. Il n'y a 
\emph{a priori} que deux points possibles pour localiser, $0$
et $\infty$, car ce sont les seuls points de $\Sr$ fix\'es par
la dilatation $z \mapsto q z$. En localisant en $0$, c'est \`a 
dire par extension de base $\Kr \hookrightarrow \Ka$, on obtient
un plongement $\Er \hookrightarrow \E$. Il y a similairement une
localisation en $\infty$: 
$\Er \hookrightarrow \mathcal{E}^{(\infty)}$ (cette derni\`ere
cat\'egorie est d'ailleurs isomorphe \`a $\E$). \\

En fait, il y aurait lieu de localiser aux \emph{singularit\'es
interm\'ediaires} de $A$, \ie\ les singularit\'es dans $\C^{*}$.
Mais elles bougent sous l'action de $z \mapsto q z$, et il faut
en fait trouver une notion de localisation en une $q$-spirale
ou en un point de $\Eq$. Le peu que nous savons faire en ce sens
est l'objet de cette section.


\subsection{Le cas des \'equations fuchsiennes}
\label{subsection:fuchsienglobal}

C'est celui o\`u Birkhoff a pos\'e et r\'esolu le probl\`eme du
recollement des donn\'ees locales. Classiquement, pour les
\'equations diff\'erentielles, c'est le probl\`eme des matrices
de connexion. La corrsepondance de Riemann-Hilbert permet, \`a
partir des monodromies locales et d'un nombre fini de matrices
de connexion, de reconstruire la classe d'\'equivalence rationnelle
d'une \'equation diff\'erentielle lin\'eaire fuchsienne. Sous forme
moderne, les donn\'ees locales et de connexion sont traduites comme
une repr\'esentation de monodromie. \\

Disons que $A \in \GL_{n}(\Kr)$ est \emph{fuchsienne} si elle l'est
en $0$ et $\infty$ (\ie\ dans $\E$ et dans $\mathcal{E}^{(\infty)}$).
Il revient au m\^eme de dire que $A$ est rationnellement \'equivalente
(\cad\ via une transformation de jauge $F \in \GL_{n}(\Kr)$ \`a une
matrice qui est non singuli\`ere en $0$ et $\infty$. Nous supposerons
donc, comme le fait Birkhoff, que $A(0), A(\infty) \in \GL_{n}(\C)$. \\

Birkhoff construit alors des ``solutions locales'' fondamentales
$\X^{(0)}$ et $\X^{(\infty)}$ de $\sq X = A X$, qui sont m\'eromorphes 
(mais multiformes) sur $\C^{*}$, et d\'efinit la \emph{matrice de
connexion de Birkhoff}:
$$
P = \left(\X^{(\infty)}\right)^{-1} \X^{(0)}.
$$
La matrice $P$ est m\'eromorphe (mais multiforme) sur $\C^{*}$
et v\'erifie $P(qz) = P(z)$ par construction (elle connecte deux
solutions d'une m\^eme \'equation, elle est donc ``$q$-constante'').
Elle est donc presque elliptique. Birkhoff d\'emontre que la donn\'ee
d'un nombre fini d'invariants locaux en $0$ et $\infty$ et de la matrice 
de connexion $P$ permet de retrouver la matrice $A$ \`a \'equivalence
rationnelle pr\`es. C'est pour prouver l'\emph{existence} d'une matrice
$A$ qu'il a invent\'e le th\'eor\`eme de factorisation dont nous avons
parl\'e au d\'ebut. \\

Le travail a \'et\'e repris dans \cite{JSAIF}, en n'utilisant que
des fonctions m\'eromorphes uniformes. Par exemple, on r\'esout
$\sq f = c f$ ($c \in \C^{*}$) \`a l'aide de 
$e_{q,c} = \theta_{q}/\theta_{q,c}$, l\`a o\`u Birkhoff utilisait
$z^{\log_{q} c}$. Ainsi, $P \in \GL_{n}\left(\M(\Eq)\right)$.
Le r\'esultat de Birkhoff prend la forme pr\'ecise d'une \'equivalence 
de cat\'egories. L'inconv\'enient de la matrice de Birkhoff (dans sa
version d'origine ou dans celle de \cite{JSAIF}) est qu'elle ne se
comporte pas bien par produit tensoriel, d'o\`u des difficult\'es
pour la th\'eorie de Galois (\cad\ pour obtenir une repr\'esentation
de groupes). Ce mauvais comportement a sa source dans le fait que
$e_{q,c} e_{q,d} \neq e_{q,cd}$ et ce, quel que soit le choix fait
des solutions m\'eromorphes $e_{q,c}$. \\

En fait, pour des \'equations \emph{r\'eguli\`eres en $0$ et $\infty$},
\cad\ telles que $A(0) = A(\infty) = I_{n}$, on n'a pas besoin de
fonctions sp\'eciales pour r\'esoudre l'\'equation, on peut le faire
avec des s\'eries convergentes. Etingof a montr\'e dans ce cas
(\cite{Etingof}) que les valeurs de la matrice de connexion engendrent
le groupe de Galois. Sans cette hypoth\`ese, van der Put et Singer ont
obtenu un r\'esultat similaire mais \`a l'aide de la r\'esolution
symbolique: les $e_{q,c}$ sont remplac\'es par des symboles $e_{c}$
tels que $e_{c} e_{d} = e_{cd}$. On n'obtient pas par cette voie de
v\'eritables invariants transcendants.

\paragraph{Avec des fibr\'es, \c{c}a marche mieux.}

Dans \cite{JSGAL}, on proc\`ede comme suit. D'apr\`es le lemme-cl\'e
de la section \ref{subsubsection:equationsfuchsiennes}, on peut \'ecrire:
\begin{eqnarray*}
A 
& = & F^{(0)}\left[A^{(0)}\right], \quad
A^{(0)} \in \GL_{n}(\C) \text{~et~} F^{(0)}(z) \in \GL_{n}(\Ka) \\
& = & F^{(\infty)}\left[A^{(\infty)}\right], \quad
A^{(\infty)} \in \GL_{n}(\C) \text{~et~} F^{(\infty)}(w) 
\in \GL_{n}(\C(\{w\})),
\end{eqnarray*}
o\`u $w = z^{-1}$ (uniformisante en $\infty \in \Sr$). Des \'equations
fonctionnelles satisfaites par $F^{(0)}$ et $F^{(\infty)}(w)$ d\'ecoule
d'ailleurs qu'ils admettent des prolongements m\'eromorphes respectivement
\`a $\C$ et \`a $\Sr \setminus \{0\}$. \\

La matrice $F =  \left(F^{(\infty)}\right)^{-1} F^{(0)}$ v\'erifie alors:
$F \in \GL_{n}(\Kw)$ et $F\left[A^{(0)}\right] = A^{(\infty)}$. Elle code
donc un \emph{isomorphisme m\'eromorphe}:
$$
\phi: \F^{(0)} = \F_{A^{(0)}} \rightarrow \F^{(\infty)} = \F_{A^{(\infty)}}.
$$
En fait, $\F^{(0)}$ est le fibr\'e local \'etudi\'e aux sections
\ref{section:constructionslocales} et \ref{section:phenomenedeStokes}.
On d\'emontre alors que $A \leadsto (\F^{(0)},\phi,\F^{(\infty)})$
est une \emph{\'equivalence de cat\'egories tannakiennes}. On en d\'eduit
que le groupe de Galois de $A$ est engendr\'e par ses composantes locales
$G_{f}^{(0)}$ et $G_{f}^{(\infty)}$, plus des donn\'ees de recollement
qui sont essentiellement les valeurs de $\phi$. Le r\'esultat d'Etingof
est un cas particulier. \\

\Rem
Un premier r\'esultat, beaucoup moins \'el\'egant, avait consist\'e
\`a tordre la matrice de connexion, de mani\`ere assez compliqu\'ee, 
pour que ses valeurs contribuent au groupe de Galois. C'est cette
version peu conceptuelle qui sert dans la pratique: pour la confluence
des automorphismes galoisiens dans \cite{JSGAL}; et pour le calcul par
Julien Roques de la plupart des groupes de Galois $q$-hyperg\'eom\'etriques 
dans \cite{Roques}. 

\paragraph{Localisation dans le cas ab\'elien.}

Les constructions ci-dessus ont toutes le m\^eme d\'efaut: il faut
une quantit\'e non-d\'enombrable de g\'en\'erateurs $P(a)$ ou $\phi(a)$
pour engendrer le groupe de Galois. Cela est \`a comparer avec la
repr\'esentation de monodromie, qui est de type fini. Le groupe
de Galois obtenu par voie alg\'ebrique n'a pas le caract\`ere discret
et transcendant de la ``vraie'' correspondance de Riemann-Hilbert. \\

Dans le cas r\'egulier, ou seule la matrice de connexion compte,
on a vu que le groupe de Galois de $A$ est param\'etr\'e par la
fonction matricielle elliptique $P$. (En fait, par la fonction
$a \mapsto P(a_{0})^{-1} P(a)$, mais on peut supposer que 
$P(a_{0}) = I_{n}$). Comme une fonction m\'eromorphe sur la
surface de Riemann $\Eq$ est la m\^eme chose qu'une fonction
rationnelle sur la courbe alg\'ebrique $\Eq$, on a donc un groupe
alg\'ebrique rationnellement param\'etr\'e par une courbe alg\'ebrique.
Dans le cas o\`u le groupe est ab\'elien, cette situation rel\`eve
de la \emph{th\'eorie g\'eom\'etrique du corps de classes}
(\cite{SerreGACC}). On a pu ainsi, dans \cite{JSGAL} localiser
la matrice de connexion et en d\'eduire un \'equivalent raisonnable
du groupe de monodromie. Cette description est trop lourde pour \^etre
reprise ici. Elle est surtout peu utile (sauf comme encouragement),
car les \'equations ab\'eliennes sont trop rares (elles sont presque
la m\^eme chose que les \'equations de rang $1$).


\subsection{Le cas g\'en\'eral}
\label{subsection:casgeneral}

La vraie g\'en\'eralisation attendue n'est pas tellement le passage
du cas fuchsien au cas irr\'egulier (qui commence \`a \^etre bien
compris), mais le passage au cas non ab\'elien: comment alors 
localiser l'effet des singularit\'es ? On va proposer une construction
int\'eressante qui r\'ealise platoniquement cette localisation. \\

Soit $A \in \GL_{n}(\Kr)$. (Une bonne partie de ce qui suit garde
un sens si $A \in \GL_{n}(\Kw)$.) Notons $Sing(A)$ le lieu singulier
de $A$, d\'efini comme sui:
$$
Sing(A) = \{\text{p\^oles de~} A \text{~dans~} \C^{*}\} \cup
\{\text{p\^oles de~} A^{-1} \text{~dans~} \C^{*}\}.
$$
Soit par ailleurs $U$ un ouvert connexe de $\C^{*}$ v\'erifiant les
deux conditions suivantes:
\begin{enumerate}
\item{$\pi(U) = \Eq$; la restriction \`a $U$ du rev\^etement
$\pi: \C^{*} \rightarrow \Eq$ est donc un isomorphisme local.
Un exemple de tel ouvert est toute couronne
$\Co(r,R) = \{z \in \C^{*} \mid r < |z| < R\}$, o\`u $R > r |q|$, en
particulier les disques \'epoint\'es $\Co(0,r)$ et $\Co(r,\infty)$
pour $r \in \R_{+}^{*}$.}
\item{$U \cap q^{-1} U \cap Sing(A) = \emptyset$; l'ouvert $U$ ne contient 
donc aucune \emph{paire singuli\`ere} $\{z,qz\} \subset Sing(A)$. (Une bonne
partie de ce qui suit reste valable sans cette condition.)}
\end{enumerate}

On pose alors:
$$
\F_{U,A} = \dfrac{U \times \C^{n}}{\sim_{A}},
$$
avec la m\^eme d\'efinition de $\sim_{A}$ que pr\'ec\'edemment.
C'est un fibr\'e vectoriel holomorphe sur $\Eq$. Le faisceau des
sections se calcule ainsi:
$$
\Gamma(V,\F_{U,A}) = 
\{X \in \O\left(U \cap \pi^{-1}(V)\right)^{n} \, \mid \, \sq X = A X\}.
$$
Une section $X \in \Gamma(V,\F_{U,A})$, vue comme fonction holomorphe
$U \cap \pi^{-1}(V)$, admet automatiquement un prolongement m\'eromorphe
\`a $\pi^{-1}(V)$, en vertu de l'\'equation fonctionnelle $\sq X = A X$.
On peut d\'ecrire le fibr\'e $\F_{U,A}$ en termes de diviseurs matriciels,
proches de \cite{Weil}. Pour cela on introduit une solution m\'eromorphe
fondamentale $\X_{0} \in \GL_{n}(\Kw)$ de $\sq X = A X$. (Il en existe
pour les m\^emes raisons qu'auparavant.) Alors le faisceau $\F_{U,A}$
est isomorphe, via la transformation de jauge $X = \X_{0} Y$, au faisceau
$\F_{U,\X_{0}}$ d\'efini par:
$$
\Gamma(V,\F_{U,\X_{0}}) = 
\{\O_{\Eq}(V)^{n} \, \mid \, 
\X_{0} Y \text{~est holomorphe sur~} U \cap \pi^{-1}(V)\}.
$$
On a not\'e ici $\O_{\Eq}$ le faisceau des fonctions holomorphes 
sur $\Eq$. Il s'identifie naturellement au faisceau
$V \mapsto \O_{\C^{*}}(\pi^{-1}(V))^{\sq}$ des fonctions holomorphes
$\sq$-invariantes sur $\C^{*}$, ce qui donne un sens au produit $\X_{0} Y$.

\paragraph{Application souhait\'ee \`a la localisation.}

Si $U$ et $U'$ sont deux ouverts v\'erifiant les conditions pr\'ec\'edentes,
il y a un isomorphisme m\'eromorphe naturel $\phi_{U,U',A}$ de $\F_{U,A}$
sur $\F_{U',A}$. Par exemple, si $U$ et $U'$ sont des disques \'epoint\'es
respectivement centr\'es en $0$ et en $\infty$, on obtient essentiellement
la matrice de connexion (sous sa forme intrins\`eque). \\

En g\'en\'eral, on peut se ramener sans difficult\'e au cas o\`u 
$Sing(A) = \{z_{1},\ldots,z_{l}\}$, avec 
$|z_{i+1}| > |q z_{i}|$ pour $1 \leq i \leq l-1$ et
$\overline{z_{i}} \neq \overline{z_{j}}$ pour $1 \leq i < j \leq l$.
On peut de plus choisir des r\'eels positifs $r_{i}$ tels que
$|z_{i}| < r_{i} < |q z_{i}|$ pour $1 \leq i \leq l$. Posons
de plus $r_{0} = 0$, $r_{l+1} = +\infty$ et $U_{i} = \Co(r_{i-1},r_{i})$
pour $1 \leq i \leq l+1$ et $\phi_{i} = \phi_{U_{i},U_{i+1},A}$ pour
$1 \leq i \leq l$. Alors chaque isomorphisme m\'eromorphe $\phi_{i}$
admet pour seule singularit\'e $\overline{z_{i}} \in \Eq$, et leur
produit est la matrice de connexion $\phi = \phi_{U_{0},U_{l+1},A}$.
On a donc localis\'e 
\footnote{On a choisi d'utiliser de vraies couronnes, \`a bords
circulaires. En fait, on pourrait prendre des couronnes topologiques,
n'imposer aucune condition aux $z_{i}$ et les parcourir dans n'importe
quel ordre. Une g\'eom\'etrie et une combinatoire de $\C^{*}$
interviennent donc ici.}
les singularit\'es de celle-ci. De plus, chaque 
$U_{i}$ permet de construire des foncteurs fibres, et les valeurs des 
$\phi_{i}$ sur $\Eq$ fournissent des op\'erateurs galoisiens. \\

Comparons maintenant le probl\`eme \`a celui du groupe de Stokes. 
On avait \'egalement une quantit\'e non d\'enombrable de g\'en\'erateurs,
les $S_{\overline{c},\overline{d}} \hat{F}_{A}(a)$. Mais ceux-ci
\'etaient tous dans un m\^eme groupe unipotent, que l'on a pu remplacer
par son alg\`ebre de Lie. Une fois le probl\`eme lin\'earis\'e
(et ab\'elianis\'e), on pouvait localiser l'effet des singulari\'es
par prise de r\'esidus. Nous ne savons rien faire de tel ici, parce que
nous ne connaissons pas de forme normale maniable pour les $\phi_{i}$. \\

Il est \`a noter que Krichever a r\'eussi dans \cite{Krichever} \`a
traiter un probl\`eme analogue pour les \'equations aux diff\'erences.



\end{document}